# CONVERGENCE RATES OF POSTERIOR DISTRIBUTIONS FOR NONIID OBSERVATIONS

By Subhashis Ghosal[1] and Aad van der Vaart

*North Carolina State University and Vrije Universiteit Amsterdam*

We consider the asymptotic behavior of posterior distributions and Bayes estimators based on observations which are required to be neither independent nor identically distributed. We give general results on the rate of convergence of the posterior measure relative to distances derived from a testing criterion. We then specialize our results to independent, nonidentically distributed observations, Markov processes, stationary Gaussian time series and the white noise model. We apply our general results to several examples of infinite-dimensional statistical models including nonparametric regression with normal errors, binary regression, Poisson regression, an interval censoring model, Whittle estimation of the spectral density of a time series and a nonlinear autoregressive model.

**1. Introduction.** Let $(\mathfrak{X}^{(n)}, \mathcal{A}^{(n)}, P_\theta^{(n)} : \theta \in \Theta)$ be a sequence of statistical experiments with observations $X^{(n)}$, where the parameter set $\Theta$ is arbitrary and $n$ is an indexing parameter, usually the sample size. We put a prior distribution $\Pi_n$ on $\theta \in \Theta$ and study the rate of convergence of the posterior distribution $\Pi_n(\cdot|X^{(n)})$ under $P_{\theta_0}^{(n)}$, where $\theta_0$ is the "true value" of the parameter. The rate of this convergence can be measured by the size of the smallest shrinking balls around $\theta_0$ that contain most of the posterior probability. For parametric models with independent and identically distributed (i.i.d.) observations, it is well known that the posterior distribution converges at the rate $n^{-1/2}$. When $\Theta$ is infinite-dimensional, but the observations are i.i.d., Ghosal, Ghosh and van der Vaart [14] obtained rates of convergence in terms of the size of the model (measured by the metric entropy or existence of certain tests) and the concentration rate of the prior around $\theta_0$ and

Received January 2004; revised March 2006.
[1]Supported in part by NSF Grant number DMS-03-49111.

computed the rate of convergence for a variety of examples. A similar result was obtained by Shen and Wasserman [27] under stronger conditions.

Little is known about the asymptotic behavior of the posterior distribution in infinite-dimensional models when the observations are not i.i.d. For independent, nonidentically distributed (i.n.i.d.) observations, consistency has recently been addressed by Amewou-Atisso, Ghosal, Ghosh and Ramamoorthi [1] and Choudhuri, Ghosal and Roy [7]. The main purpose of the present paper is to obtain a theorem on rates of convergence of posterior distributions in a general framework not restricted to the setup of i.i.d. observations. We specialize this theorem to several classes of non-i.i.d. models including i.n.i.d. observations, Gaussian time series, Markov processes and the white noise model. The theorem applies in every situation where it is possible to test the true parameter versus balls of alternatives with exponential error probabilities and it is not restricted to any particular structure on the joint distribution. The existence of such tests has been proven in many special cases by Le Cam [20, 21, 22] and Birgé [3, 4, 5], who used them to construct estimators with optimal rates of convergence, determined by the (local) metric entropy or "Le Cam dimension" of the model. Our main theorem uses the same metric entropy measure of the complexity of the model and combines this with a measure of prior concentration around the true parameter to obtain a bound on the posterior rate of convergence, generalizing the corresponding result of Ghosal, Ghosh and van der Vaart [14]. We apply these results to obtain posterior convergence rates for linear regression, nonparametric regression, binary regression, Poisson regression, interval censoring, spectral density estimation and nonlinear autoregression. van der Meulen, van der Vaart and van Zanten [30] have extended the approach of this paper to several types of diffusion models.

The organization of the paper is as follows. In the next section, we describe our main theorem in an abstract framework. In Sections 3, 4, 5 and 6, we specialize to i.n.i.d. observations, Markov chains, the white noise model and Gaussian time series, respectively. In Section 7, we discuss a large number of more concrete applications, combining models of various types with many types of different priors, including priors based on the Dirichlet process, mixture representations or sequence expansions on spline bases, priors supported on finite sieves and conjugate Gaussian priors. Technical proofs, including the proofs of the main results, are collected in Section 8.

The notation $\lesssim$ will be used to denote inequality up to a constant that is fixed throughout. The notation $Pf$ will abbreviate $\int f \, dP$. The symbol $\lfloor x \rfloor$ will stand for the greatest integer less than or equal to $x$. Let $h(f,g) = (\int (f^{1/2} - g^{1/2})^2 \, d\mu)^{1/2}$ and $K(f,g) = \int f \log(f/g) \, d\mu$ stand for the Hellinger distance and Kullback–Leibler divergence, respectively, between two nonnegative densities $f$ and $g$ relative to a measure $\mu$. Furthermore, we define additional discrepancy measures by $V_k(f,g) = \int f |\log(f/g)|^k \, d\mu$ and $V_{k,0}(f,g) = \int f|\log(f/g) - K(f,g)|^k \, d\mu$, $k > 1$. The index $k = 2$ of $V_2$ and $V_{2,0}$ may be omitted and these simply written as $V$ and $V_0$, respectively. The symbols $\mathbb{N}$ and $\mathbb{R}$ will denote the sets of natural and real numbers, respectively. The $\varepsilon$-covering number of a set $\Theta$ for a semimetric $d$, denoted by $N(\varepsilon, \Theta, d)$, is the minimal number of $d$-balls of radius $\varepsilon$ needed to cover the set $\Theta$; see, for example, [31].

**2. General theorem.** For each $n \in \mathbb{N}$ and $\theta \in \Theta$, let $P_\theta^{(n)}$ admit densities $p_\theta^{(n)}$ relative to a $\sigma$-finite measure $\mu^{(n)}$. Assume that $(x,\theta) \mapsto p_\theta^{(n)}(x)$ is jointly measurable relative to $\mathcal{A} \otimes \mathcal{B}$, where $\mathcal{B}$ is a $\sigma$-field on $\Theta$. By Bayes' theorem, the posterior distribution is given by

$$(2.1) \qquad \Pi_n(B|X^{(n)}) = \frac{\int_B p_\theta^{(n)}(X^{(n)}) \, d\Pi_n(\theta)}{\int_\Theta p_\theta^{(n)}(X^{(n)}) \, d\Pi_n(\theta)}, \qquad B \in \mathcal{B}.$$

Here, $X^{(n)}$ is an "observation," which, in our setup, will be understood to be generated according to $P_{\theta_0}^{(n)}$ for some given $\theta_0 \in \Theta$.

For each $n$, let $d_n$ and $e_n$ be semimetrics on $\Theta$ with the property that there exist universal constants $\xi > 0$ and $K > 0$ such that for every $\varepsilon > 0$ and for each $\theta_1 \in \Theta$ with $d_n(\theta_1, \theta_0) > \varepsilon$, there exists a test $\phi_n$ such that

$$(2.2) \qquad P_{\theta_0}^{(n)} \phi_n \leq e^{-Kn\varepsilon^2}, \qquad \sup_{\theta \in \Theta: e_n(\theta, \theta_1) < \varepsilon \xi} P_\theta^{(n)}(1 - \phi_n) \leq e^{-Kn\varepsilon^2}.$$

Typically, we have $d_n \leq e_n$ and in many cases we choose $d_n = e_n$, but using two semimetrics provides some added flexibility. Le Cam [20, 21, 22] and Birgé [3, 4, 5] showed that the rate of convergence, in a minimax sense, of the best estimators of $\theta$ relative to the distance $d_n$ can be understood in terms of the *Le Cam dimension* or *local entropy* function of the set $\Theta$ relative to $d_n$. For our purposes, this dimension is a function whose value at $\varepsilon > 0$ is defined to be $\log N(\varepsilon\xi, \{\theta : d_n(\theta, \theta_0) \leq \varepsilon\}, e_n)$, that is, the logarithm of the minimum number of $d_n$-balls of radius $\varepsilon\xi$ needed to cover an $e_n$-ball of radius $\varepsilon$ around the true parameter $\theta_0$. Birgé [3, 4] and Le Cam [20, 21, 22] showed that there exist estimators $\hat{\theta}_n = \hat{\theta}_n(X^{(n)})$ such that $d_n(\hat{\theta}_n, \theta_0) = O_P(\varepsilon_n)$ under $P_{\theta_0}^{(n)}$, where

$$(2.3) \qquad \sup_{\varepsilon > \varepsilon_n} \log N(\varepsilon\xi, \{\theta : d_n(\theta, \theta_0) \leq \varepsilon\}, e_n) \leq n\varepsilon_n^2.$$

Further, under certain conditions $\varepsilon_n$ is the best rate obtainable, given the model, and hence gives a minimax rate.

As in the i.i.d. case, the behavior of posterior distributions depends on the size of the model measured by (2.3) and the concentration rate of the prior $\Pi_n$ at $\theta_0$. For a given $k > 1$, let

$$B_n(\theta_0, \varepsilon; k) = \{\theta \in \Theta : K(p_{\theta_0}^{(n)}, p_\theta^{(n)}) \leq n\varepsilon^2, V_{k,0}(p_{\theta_0}^{(n)}, p_\theta^{(n)}) \leq n^{k/2}\varepsilon^k\}.$$

An appropriate condition will appear as a lower bound on $\Pi_n(B_n(\theta_0;\varepsilon,k))$ with $k=2$ being good enough to establish convergence in mean. For almost sure convergence, or convergence of the posterior mean, better control may be needed (through a larger value of $k$), depending on the rate of convergence.

The following result, generalizing Theorem 2.4 of Ghosal, Ghosh and van der Vaart [14] for the i.i.d. case, bounds the rate of posterior convergence.

THEOREM 1. *Let $d_n$ and $e_n$ be semimetrics on $\Theta$ for which tests satisfying* (2.2) *exist. Let $\varepsilon_n > 0$, $\varepsilon_n \to 0$, $(n\varepsilon_n^2)^{-1} = O(1)$, $k > 1$, and $\Theta_n \subset \Theta$ be such that for every sufficiently large $j \in \mathbb{N}$,*

$$(2.4) \qquad \sup_{\varepsilon > \varepsilon_n} \log N\left(\frac{1}{2}\varepsilon\xi, \{\theta \in \Theta_n : d_n(\theta,\theta_0) < \varepsilon\}, e_n\right) \leq n\varepsilon_n^2,$$

$$(2.5) \qquad \frac{\Pi_n(\theta \in \Theta_n : j\varepsilon_n < d_n(\theta,\theta_0) \leq 2j\varepsilon_n)}{\Pi_n(B_n(\theta_0,\varepsilon_n;k))} \leq e^{Kn\varepsilon_n^2 j^2/2}.$$

*Then for every $M_n \to \infty$, we have that*

$$(2.6) \qquad P_{\theta_0}^{(n)}\Pi_n(\theta \in \Theta_n : d_n(\theta,\theta_0) \geq M_n\varepsilon_n|X^{(n)}) \to 0.$$

The theorem uses the fact that $\Theta_n \subset \Theta$ to alleviate the entropy condition (2.4), but returns an assertion about the posterior distribution on $\Theta_n$ only. The complementary assertion $P_{\theta_0}^{(n)}\Pi_n(\Theta \setminus \Theta_n|X^{(n)}) \to 0$ may be handled either by a direct argument or by the following analog of Lemma 5 of [2].

LEMMA 1. *If $\frac{\Pi_n(\Theta\setminus\Theta_n)}{\Pi_n(B_n(\theta_0,\varepsilon_n;k))} = o(e^{-2n\varepsilon_n^2})$ for some $k > 1$, then $P_{\theta_0}^{(n)}\Pi_n(\Theta \setminus \Theta_n|X^{(n)}) \to 0$.*

The choice $\Theta_n = \Theta$, which makes the condition of Lemma 1 trivial, imposes a much stronger restriction on (2.4) and is generally unattainable when $\Theta$ is not compact.

The following theorem extends the convergence in Theorem 1 to almost sure convergence and yields a rate for the convergence under slightly stronger conditions.

THEOREM 2. *In the situation of Theorem 1,*

(i) *if all $X^{(n)}$ are defined on a fixed sample space and $\varepsilon_n \gtrsim n^{-\alpha}$ for some $\alpha \in (0,1/2)$ such that $k(1-2\alpha) > 2$, then the convergence* (2.6) *also holds in the almost sure sense;*

(ii) *if $\varepsilon_n \gtrsim n^{-\alpha}$ for some $\alpha \in (0,1/2)$ such that $k(1-2\alpha) > 4\alpha$, then the left side of* (2.6) *is $O(\varepsilon_n^2)$.*

If $\Theta$ is a convex set and $d_n^2$ is a convex function in one argument keeping the other fixed and is bounded above by $B$, then for $\hat{\theta}_n = \int \theta \, d\Pi_n(\theta|X^{(n)})$, we have, by Jensen's inequality, that

$$d_n^2(\hat{\theta}_n, \theta_0) \leq \int d_n^2(\theta, \theta_0) \, d\Pi_n(\theta|X^{(n)}) \leq \varepsilon_n^2 + B^2 \Pi_n(d_n(\theta, \theta_0) \geq \varepsilon_n|X^{(n)}).$$

This yields the rate $\varepsilon_n$ for the point estimator $\hat{\theta}_n$ under the conditions of Theorem 1.

The complicated-looking condition (2.5) can often be simplified in infinite-dimensional cases, where, typically, $n\varepsilon_n^2 \to \infty$. Because the numerator in (2.5) is trivially bounded by one, a sufficient condition for (2.5) is that $\Pi_n(B_n(\theta_0, \varepsilon_n, k)) \gtrsim e^{-cn\varepsilon_n^2}$. The local entropy in condition (2.4) can also often be replaced by the global entropy $\log N(\varepsilon\xi/2, \Theta_n, e_n)$ without affecting rates. Also, if the prior is such that the minimax rate given by (2.3) satisfies (2.5) and the condition of Lemma 1, then the posterior convergence rate attains the minimax rate.

Entropy conditions, however, may not always be appropriate to ensure the existence of tests. Ad hoc tests may sometimes be more conveniently constructed. A more general theorem on convergence rates, which is formulated directly in terms of tests and stated below, may be proven in a similar manner.

THEOREM 3. *Let $d_n$ be a semimetric on $\Theta$, $\varepsilon_n \to 0$, $(n\varepsilon_n^2)^{-1} = O(1)$, $k > 1$, $K > 0$, $\Theta_n \subset \Theta$ and $\phi_n$ be a sequence of test functions such that*

$$(2.7) \quad P_{\theta_0}^{(n)} \phi_n \to 0, \qquad \sup_{\theta \in \Theta_n : j\varepsilon_n < d_n(\theta, \theta_0) \leq 2j\varepsilon_n} P_\theta^{(n)}(1 - \phi_n) \lesssim e^{-Kj^2 n\varepsilon_n^2}$$

*and (2.5) holds. Then for every $M_n \to \infty$, we have that $P_{\theta_0}^{(n)} \Pi_n(\theta \in \Theta_n : d_n(\theta, \theta_0) \geq M_n \varepsilon_n | X^{(n)}) \to 0$.*

**3. Independent observations.** In this section, we consider the case where the observation $X^{(n)}$ is a vector $X^{(n)} = (X_1, X_2, \ldots, X_n)$ of independent observations $X_i$. Thus, we take the measures $P_\theta^{(n)}$ of Section 2 equal to product measures $\bigotimes_{i=1}^n P_{\theta,i}$ on a product measurable space $\bigotimes_{i=1}^n (\mathfrak{X}_i, \mathcal{A}_i)$. We assume that the distribution $P_{\theta,i}$ of the $i$th component $X_i$ possesses a density $p_{\theta,i}$ relative to a $\sigma$-finite measure $\mu_i$ on $(\mathfrak{X}_i, \mathcal{A}_i)$, $i = 1, \ldots, n$. In this case, tests can be constructed relative to the semimetric $d_n$, whose square is given by

$$(3.1) \qquad d_n^2(\theta, \theta') = \frac{1}{n} \sum_{i=1}^n \int (\sqrt{p_{\theta,i}} - \sqrt{p_{\theta',i}})^2 \, d\mu_i.$$

Thus, $d_n^2$ is the average of the squares of the Hellinger distances for the distributions of the individual observations.

The following lemma, due to Birgé (cf. [22], page 491, or [4], Corollary 2 on page 149), guarantees the existence of tests satisfying the conditions of (2.2).

LEMMA 2. *If $P_\theta^{(n)}$ are product measures and $d_n$ is defined by* (3.1), *then there exist tests $\phi_n$ such that $P_{\theta_0}^{(n)}\phi_n \leq e^{-\frac{1}{2}nd_n^2(\theta_0,\theta_1)}$ and $P_\theta^{(n)}(1-\phi_n) \leq e^{-\frac{1}{2}nd_n^2(\theta_0,\theta_1)}$ for all $\theta \in \Theta$ such that $d_n(\theta,\theta_1) \leq \frac{1}{18}d_n(\theta_0,\theta_1)$.*

The Kullback–Leibler divergence between product measures is equal to the sum of the Kullback–Leibler divergences between the individual components. Furthermore, as a consequence of the Marcinkiewiz–Zygmund inequality (e.g., [9], page 356), the mean $\bar{Y}_n$ of $n$ independent random variables satisfies $E|\bar{Y}_n - E\bar{Y}_n|^k \leq C_k n^{-k/2} \frac{1}{n}\sum_{i=1}^n E|Y_i|^k$ for $k \geq 2$, where $C_k$ is a constant depending only on $k$. Therefore, the set $B_n(\theta_0, \varepsilon; k)$ contains the set

$$B_n^*(\theta_0, \varepsilon; k) = \left\{\theta \in \Theta : \frac{1}{n}\sum_{i=1}^n K_i(\theta_0, \theta) \leq \varepsilon^2, \ \frac{1}{n}\sum_{i=1}^n V_{k,0;i}(\theta_0, \theta) \leq C_k \varepsilon^k\right\},$$

where $K_i(\theta_0, \theta) = K(P_{\theta_0,i}, P_{\theta,i})$ and $V_{k,0;i}(\theta_0, \theta) = V_{k,0}(P_{\theta_0,i}, P_{\theta,i})$. Thus, we can work with a "ball" around $\theta_0$ relative to the average Kullback–Leibler divergence and the average $k$th order moments, as in the preceding display, and simplify Theorem 1 to the following result:

THEOREM 4. *Let $P_\theta^{(n)}$ be product measures and $d_n$ be defined by* (3.1). *Suppose that for a sequence $\varepsilon_n \to 0$ such that $n\varepsilon_n^2$ is bounded away from zero, some $k > 1$, all sufficiently large $j$ and sets $\Theta_n \subset \Theta$, the following conditions hold:*

(3.2) $$\sup_{\varepsilon > \varepsilon_n} \log N(\varepsilon/36, \{\theta \in \Theta_n : d_n(\theta, \theta_0) < \varepsilon\}, d_n) \leq n\varepsilon_n^2;$$

(3.3) $$\frac{\Pi_n(\Theta \setminus \Theta_n)}{\Pi_n(B_n^*(\theta_0, \varepsilon_n; k))} = o(e^{-2n\varepsilon_n^2});$$

(3.4) $$\frac{\Pi_n(\theta \in \Theta_n : j\varepsilon_n < d_n(\theta, \theta_0) \leq 2j\varepsilon_n)}{\Pi_n(B_n^*(\theta_0, \varepsilon_n; k))} \leq e^{n\varepsilon_n^2 j^2/4}.$$

*Then $P_{\theta_0}^{(n)}\Pi_n(\theta : d_n(\theta, \theta_0) \geq M_n\varepsilon_n|X^{(n)}) \to 0$ for every $M_n \to \infty$.*

The average Hellinger distance is not always the most natural choice. It can be replaced by any other distance $d_n$ that satisfies (3.2)–(3.3) and for which the conclusion of Lemma 2 holds. Often, we set $k=2$ and work with the smaller neighborhood

$$(3.5) \quad \bar{B}_n(\theta_0, \varepsilon) = \left\{\theta : \frac{1}{n}\sum_{i=1}^n K_i(\theta_0, \theta) \leq \varepsilon^2, \frac{1}{n}\sum_{i=1}^n V_{2;i}(\theta_0, \theta) \leq \varepsilon^2\right\}.$$

**4. Markov chains.** For $\theta$ ranging over a set $\Theta$, let $(x,y) \mapsto p_\theta(y|x)$ be a collection of transition densities from a measurable space $(\mathfrak{X}, \mathcal{A})$ into itself, relative to some reference measure $\nu$. Thus, for each $\theta \in \Theta$, the map $(x,y) \mapsto p_\theta(y|x)$ is measurable and for each $x$, the map $y \mapsto p_\theta(y|x)$ is a probability density relative to $\mu$. Let $X_0, X_1, \ldots$ be a stationary Markov chain generated according to the transition density $p_\theta$, where it is assumed that there exists a stationary distribution $Q_\theta$ with $\mu$-density $q_\theta$. Let $P_\theta^{(n)}$ be the law of $(X_0, X_1, \ldots, X_n)$.

Tests satisfying the conditions of (2.2) can be obtained from results of Birgé [4], which are more refined versions of his own results in [3]. A special case is presented as Lemma 3 below. Actually, Birgé's result ([4], Theorem 3, page 155) is much more general in that it also applies to nonstationary chains and allows different upper and lower bounds, as seen in the following display.

Assume that there exists a finite measure $\nu$ on $(\mathfrak{X}, \mathcal{A})$ such that, for some $k, l \in \mathbb{N}$, every $\theta \in \Theta$ and every $x \in \mathfrak{X}$ and $A \in \mathcal{A}$,

$$(4.1) \quad P_\theta(X_l \in A | X_0 = x) \lesssim \nu(A) \lesssim \frac{1}{k}\sum_{j=1}^k P_\theta(X_j \in A | X_0 = x),$$

where $P_\theta$ is the generic notation for any probability law governed by $\theta$. For instance, if there exists a $\mu$-integrable function $r$ such that $r(y) \lesssim p_\theta(y|x) \lesssim r(y)$ for every $(x,y)$, then (4.1) holds with the measure $\nu$ given by $d\nu(y) = r(y)\, d\mu(y)$. Define the square of a semidistance $d$ by

$$(4.2) \quad d^2(\theta, \theta') = \iint \left[\sqrt{p_\theta(y|x)} - \sqrt{p_{\theta'}(y|x)}\right]^2 d\mu(y)\, d\nu(x).$$

LEMMA 3. *If there exist $k$, $l$ and a measure $\nu$ such that (4.1) holds, then there exist a constant $K$ depending only on $(k,l)$ and tests $\phi_n$ such that*

$$P_{\theta_0}^{(n)} \phi_n \leq e^{-Knd^2(\theta_0, \theta_1)}, \qquad \sup_{\theta \in \Theta : d(\theta, \theta_1) \leq d(\theta_0, \theta_1)/8} P_\theta^{(n)}(1 - \phi_n) \leq e^{-Knd^2(\theta_0, \theta_1)}.$$

The preceding lemma is also true if the chain is not started at stationarity. If, as we assume, $X_0$ is generated from a stationary distribution under $\theta_0$, then the Kullback–Leibler divergence of $P_{\theta_0}^{(n)}$ and $P_\theta^{(n)}$ satisfies

$$(4.3) \quad K(P_{\theta_0}^{(n)}, P_\theta^{(n)}) = n\int K(p_{\theta_0}(\cdot|x), p_\theta(\cdot|x))\, dQ_{\theta_0}(x) + K(q_{\theta_0}, q_\theta).$$

To handle the neighborhoods $B_n(\theta_0, \varepsilon; 2)$, we need a bound on $V(P_{\theta_0}^{(n)}, P_\theta^{(n)})$, which will also be of the order of $n$ times an expression depending only on individual observations, under a variety of conditions. In the following lemma, we use an $\alpha$-mixing assumption. For a sequence $\{X_n\}$, let the $\alpha$-mixing coefficient be given by $\alpha_h = \sup\{|\Pr(X_0 \in A, X_h \in B) - \Pr(X_0 \in A)\Pr(X_h \in B)| : A, B \in \mathcal{B}(\mathbb{R})\}$.

LEMMA 4. *Suppose that the Markov chain $X_0, X_1, \ldots$ is $\alpha$-mixing under $\theta_0$, with mixing coefficients $\alpha_h$. Then for every $s > 2$, $V(p_{\theta_0}^{(n)}, p_\theta^{(n)})$ is bounded by*

$$\frac{8sn}{s-2} \sum_{h=0}^\infty \alpha_h^{1-2/s} \left( \iint \left| \log \frac{p_{\theta_0}(y|x)}{p_\theta(y|x)} \right|^s p_{\theta_0}(y|x)\, d\mu(y)\, dQ_{\theta_0}(x) \right)^{2/s} + 2V(q_{\theta_0}, q_\theta).$$

PROOF. We can write

$$(4.4) \qquad \log \frac{p_{\theta_0}^{(n)}}{p_\theta^{(n)}} = \sum_{i=1}^n \log \frac{p_{\theta_0}(X_i|X_{i-1})}{p_\theta(X_i|X_{i-1})} + \log \frac{q_{\theta_0}(X_0)}{q_\theta(X_0)} =: n\bar{Y}_n + Z_0,$$

where $Y_i = \log(p_{\theta_0}(X_i|X_{i-1})/p_\theta(X_i|X_{i-1}))$ and $Z_0 = \log(q_{\theta_0}(X_0)/q_\theta(X_0))$. Then $Y_1, Y_2, \ldots$ are $\alpha$-mixing with mixing coefficients $\alpha_{h-1}$. Therefore, the variance of the left-hand side of (4.4) is bounded above by $n(\mathrm{E}|Y_i|^s)^{2/s} \times 4s(s-2)^{-1} \sum_{h=1}^\infty \alpha_{h-1}^{1-2/s}$, by the bound of Ibragimov [18]. □

Let $\Theta_1 \subset \Theta$ be the set of parameter values such that $K(q_{\theta_0}, q_\theta)$ and $V(q_{\theta_0}, q_\theta)$ are bounded by 1. Then from (4.3) and Lemma 4, it follows that for large $n$ and $\varepsilon^2 \geq 2/n$, the set $B_n(\theta_0, \varepsilon; 2)$ contains the set $B^*(\theta_0, \varepsilon; s)$ defined by

$$\left\{ \theta \in \Theta_1 : \mathrm{P}_{\theta_0} \log \left( \frac{p_{\theta_0}}{p_\theta}(X_1|X_0) \right) \leq \frac{1}{2}\varepsilon^2, \mathrm{P}_{\theta_0} \left| \log \frac{p_{\theta_0}}{p_\theta}(X_1|X_0) \right|^s \leq C_s \varepsilon^s \right\},$$

where the power $s$ must be chosen sufficiently large to ensure that the mixing coefficients satisfy $\sum_{h=0}^\infty \alpha_h^{1-2/s} < \infty$ and where $C_s^{-2/s} = 16s(2-s)^{-1} \sum_{h=0}^\infty \alpha_h^{1-2/s}$. The contributions of $Q_{\theta_0}(\log(q_{\theta_0}/q_\theta))$ and $Q_{\theta_0}(\log(q_{\theta_0}/q_\theta))^2$ may also be incorporated into the bound.

The above facts may be combined to obtain the following result.

THEOREM 5. *Let $P_\theta^{(n)}$ be the distribution of $(X_0, X_1, \ldots, X_n)$ for a stationary Markov chain $X_0, X_1, \ldots$ with transition densities $p_\theta(y|x)$ and stationary density $q_\theta$ satisfying (4.1) and let $d$ be defined by (4.2). Assume, further, that the chain is $\alpha$-mixing with coefficients $\alpha_h$ satisfying $\sum_{h=0}^\infty \alpha_h^{1-1/s} < \infty$ for some $s > 2$. Suppose that for a sequence $\varepsilon_n \to 0$ such that $n\varepsilon_n^2 \geq 2$,* *some* $s > 2$, *every sufficiently large* $j$ *and sets* $\Theta_n \subset \Theta$, *the following conditions are satisfied:*

$$\sup_{\varepsilon > \varepsilon_n} \log N(\varepsilon/16, \{\theta \in \Theta_n : d(\theta, \theta_0) < \varepsilon\}, d) \leq n\varepsilon_n^2; \tag{4.5}$$

$$\frac{\Pi_n(\Theta \setminus \Theta_n)}{\Pi_n(B^*(\theta_0, \varepsilon_n; s))} = o(e^{-2n\varepsilon_n^2}); \tag{4.6}$$

$$\frac{\Pi_n(\theta \in \Theta_n : (j-1)\varepsilon_n < d(\theta, \theta_0) \leq j\varepsilon_n)}{\Pi_n(B^*(\theta_0, \varepsilon_n; s))} \leq e^{Kn\varepsilon_n^2 j^2/8}, \tag{4.7}$$

*for the constant* $K$ *of Lemma* 3. *Then* $P_{\theta_0}^{(n)} \Pi_n(\theta : d_*(\theta, \theta_0) \geq M_n \varepsilon_n | X^{(n)}) \to 0$ *for every* $M_n \to \infty$.

A Markov chain with $n$-step transition probability $P^n(x, \cdot) = \Pr(X_n \in A | X_0 = x)$ and stationary measure $Q$ is called *uniformly ergodic* if $\|P^n(x, \cdot) - Q\| \to 0$ as $n \to \infty$, uniformly in $x$, where $\|\cdot\|$ is the total variation norm. It can be shown that the convergence is then automatically exponentially fast (cf. [23], Theorem 16.0.2). Thus, the $\alpha$-mixing coefficients are exponentially decreasing and hence satisfy $\sum_{h=0}^\infty \alpha_h^{1-2/s} < \infty$ for every $s > 2$. Hence, it suffices to verify (4.7) with some arbitrary fixed $s > 2$. If $\sup\{\int |p_{\theta_0}(y|x_1) - p_{\theta_0}(y|x_2)| \, d\mu(y) : x_1, x_2 \in \mathbb{R}\} < 2$, then integrating out $x_2$ relative to the stationary measure $q_{\theta_0}$, we see that Condition (16.8) of Theorem 16.0.2 of [23] holds and hence the chain is uniformly ergodic.

**5. White noise model.** Let $\Theta \subset L_2[0, 1]$ and for $\theta \in \Theta$, let $P_\theta^{(n)}$ be the distribution on $C[0, 1]$ of the stochastic process $X^{(n)} = (X_t^{(n)} : 0 \leq t \leq 1)$ defined structurally as $X_t^{(n)} = \int_0^t \theta(s) \, ds + \frac{1}{\sqrt{n}} W_t$ for a standard Brownian motion $W$. This is the standard white noise model, which is known to arise as an approximation of many particular sequences of experiments. An equivalent experiment is obtained by the one-to-one correspondence of $X^{(n)}$ with the sequence defined by $X_{n,i} = \langle X^{(n)}, e_i \rangle$, where $\langle \cdot, \cdot \rangle$ is the inner product of $L_2[0, 1]$ and $\{e_1, e_2, \ldots\}$ is a given orthonormal basis of $L_2[0, 1]$. The variables $X_{n,1}, X_{n,2}, \ldots$ are independent and normally distributed, with means $\langle \theta, e_i \rangle$ and variance $n^{-1}$. In the following, we use this concrete representation and abuse notation by identifying $X^{(n)}$ with the sequence $(X_{n,1}, X_{n,2}, \ldots)$ and $\theta \in \Theta$ with the sequence $(\theta_1, \theta_2, \ldots)$ defined by $\theta_i = \langle \theta, e_i \rangle$. In the latter representation, we have that $\Theta \subset \ell_2$, the space of square summable sequences. Let $\|\theta\|^2 = \int_0^1 \theta^2(s) \, ds = \sum_{i=1}^\infty \theta_i^2$ denote the squared $L_2$-norm.

Tests satisfying the conditions of (2.2) can easily be found explicitly, namely, as the likelihood ratio test for $\theta_0$ versus $\theta_1$, where we can use the $L_2$-norm for both $d_n$ and $e_n$. Furthermore, the Kullback–Leibler divergence and discrepancy $V_{2,0}$ also turn out to be multiples of the $L_2$-norm.

LEMMA 5. *The test $\phi_n = 1\{2\langle \theta_1 - \theta_0, X^{(n)}\rangle > \|\theta_1\|^2 - \|\theta_0\|^2\}$ satisfies $P_{\theta_0}^{(n)} \phi_n \leq 1 - \Phi(\sqrt{n}\|\theta_1 - \theta_0\|/2)$ and $P_\theta^{(n)}(1 - \phi_n) \leq 1 - \Phi(\sqrt{n}\|\theta_1 - \theta_0\|/4)$ for any $\theta \in \Theta$ such that $\|\theta - \theta_1\| \leq \|\theta_1 - \theta_0\|/4$.*

LEMMA 6. *For every $\theta, \theta_0 \in \Theta \subset L_2[0,1]$, we have $K(P_{\theta_0}^{(n)}, P_\theta^{(n)}) = \frac{1}{2}n\|\theta - \theta_0\|^2$ and $V_{2,0}(P_{\theta_0}^{(n)}, P_\theta^{(n)}) = n\|\theta - \theta_0\|^2$. Consequently, we have $B_n(\theta_0, \varepsilon; 2) = \{\theta \in \Theta : \|\theta - \theta_0\| \leq \varepsilon\}$.*

PROOF OF LEMMA 5. The test rejects the null hypothesis for positive values of the statistic $T_n = \langle \theta_1 - \theta_0, X^{(n)}\rangle - \frac{1}{2}\|\theta_1\|^2 + \frac{1}{2}\|\theta_0\|^2$, which, under $\theta$, is distributed as $\langle \theta_1 - \theta_0, \theta - \theta_1\rangle + \frac{1}{2}\|\theta_1 - \theta_0\|^2 + \frac{1}{\sqrt{n}}\langle \theta_1 - \theta_0, W\rangle$. The variable $\langle \theta_1 - \theta_0, W\rangle$ is normally distributed with mean zero and variance $\|\theta_1 - \theta_0\|^2$. Under $\theta = \theta_0$, the mean of the test statistic is equal to $-\frac{1}{2}\|\theta_0 - \theta_1\|^2$, whereas for $\|\theta - \theta_1\| \leq \xi\|\theta_1 - \theta_0\|$ and $\xi \in (0, \frac{1}{2})$, the mean of the statistic under $\theta$ is bounded below by $(\frac{1}{2} - \xi)\|\theta_0 - \theta_1\|^2$, in view of the Cauchy–Schwarz inequality. The lemma follows upon choosing $\xi = 1/4$. □

PROOF OF LEMMA 6. We write $\log(p_{\theta_0}^{(n)}/p_\theta^{(n)}) = n\langle \theta_0 - \theta, X^{(n)}\rangle - \frac{n}{2}\|\theta_0\|^2 + \frac{n}{2}\|\theta\|^2$, whence the mean and variance under $\theta_0$ are easily obtained. □

In the preceding lemmas, no restriction on the parameter set $\Theta \subset L_2[0,1]$ was imposed. The lemmas lead to the following theorem, which gives bounds on the rate of convergence in terms of quantities involving the $L_2$-norm only.

THEOREM 6. *Let $P_\theta^{(n)}$ be the distribution on $C[0,1]$ of the solution of the diffusion equation $dX_t = \theta(t)\,dt + n^{-1/2}\,dW_t$ with $X_0 = 0$. Suppose that for $\varepsilon_n \to 0$, $(n\varepsilon_n^2)^{-1} = O(1)$ and $\Theta \subset L_2[0,1]$, the following conditions are satisfied:*

(5.1) $$\sup_{\varepsilon > \varepsilon_n} \log N(\varepsilon/8, \{\theta \in \Theta : \|\theta - \theta_0\| < \varepsilon\}, \|\cdot\|) \leq n\varepsilon_n^2;$$

*for every $j \in \mathbb{N}$*

(5.2) $$\frac{\Pi_n(\theta \in \Theta : \|\theta - \theta_0\| \leq j\varepsilon_n)}{\Pi_n(\theta \in \Theta : \|\theta - \theta_0\| \leq \varepsilon_n)} \leq e^{n\varepsilon_n^2 j^2/64}.$$

*Then $P_{\theta_0}^{(n)}\Pi_n(\theta \in \Theta : \|\theta - \theta_0\| \geq M_n\varepsilon_n | X^{(n)}) \to 0$ for every $M_n \to \infty$.*

In Section 7.6, we shall calculate the rate of convergence for a conjugate prior.

**6. Gaussian time series.** Suppose that $X_1, X_2, \ldots$ is a stationary Gaussian process with mean zero and spectral density $f$, which is known to belong to a model $\mathcal{F}$. Let $\gamma_h(f) = \int_{-\pi}^{\pi} e^{ih\lambda} f(\lambda) \, d\lambda$ be the corresponding autocovariance function. Let $P_f^{(n)}$ be the distribution of $(X_1, \ldots, X_n)$.

For this situation, we can derive the following lemma from [3]. Let $\|f\|_2$ and $\|f\|_\infty$ be the $L_2$-norm relative to Lebesgue measure and the uniform norm of a function $f : (-\pi, \pi] \to \mathbb{R}$, respectively.

LEMMA 7. *Suppose that there exist constants $\Gamma$ and $M$ such that $\|\log f\|_\infty \leq \Gamma$ and $\sum_{h=-\infty}^{\infty} |h| \gamma_h^2(f) \leq M$ for every $f \in \mathcal{F}$. Then there exist constants $\xi$ and $K$ depending only on $\Gamma$ and $M$ such that for every $\varepsilon \gtrsim 1/\sqrt{n}$ and every $f_0, f_1 \in \mathcal{F}$ with $\|f_1 - f_0\|_2 \geq \xi\varepsilon$, we have*

$$(6.1) \qquad P_{f_0}^{(n)} \phi_n \vee \sup_{f \in \mathcal{F}: \|f - f_1\|_\infty \leq \xi\varepsilon} P_f^{(n)}(1 - \phi_n) \leq e^{-Kn\varepsilon^2}.$$

PROOF. It follows from the assumptions that $\sum_{|h|>n/2} \gamma_h^2(f) \leq 2M/n$. This is bounded by $\varepsilon^2$ for $\varepsilon \geq \sqrt{2M/n}$. The assertion follows from Proposition 5.5, page 222 of [3], with $\phi_n = 1\{\log(p_{f_1}^{(n)}/p_{f_0}^{(n)}) \geq 0\}$. □

The preceding lemma shows that tests satisfying the conditions of (2.2) exist when $d_n$ is the $L_2$-distance and when $e_n$ is the uniform distance, leading to conditions in terms of $N(\varepsilon\xi, \{f \in \mathcal{F} : \|f - f_0\|_2 < \varepsilon\}, \|\cdot\|_\infty)$. We do not know if the $L_\infty$-distance can be replaced by the $L_2$-distance. The uniform bound on $\|\log f\|_\infty$ is not unreasonable as it is known that the structure of the time series changes dramatically if the spectral density approaches zero. The following lemma allows the neighborhoods $B_n(f_0, \varepsilon; 2)$ to be dealt with entirely in terms of balls for the $L_2$-norm.

LEMMA 8. *Suppose that there exists constant $\Gamma$ such that $\|\log f\|_\infty \leq \Gamma$ for every $f \in \mathcal{F}$. Then there exists a constant $C$ depending only on $\Gamma$ such that for every $f, g \in \mathcal{F}$, we have $P_f^{(n)}(\log(p_f^{(n)}/p_g^{(n)})) \lesssim Cn\|f - g\|_2^2$ and $\mathrm{var}_{P_f^{(n)}}(\log(p_f^{(n)}/p_g^{(n)})) \lesssim Cn\|f - g\|_2^2$.*

PROOF. The $(k, l)$th element of the covariance matrix $T_n(f)$ of $X^{(n)} = (X_1, \ldots, X_n)$, given the spectral density $f$, is given by $\int_{-\pi}^{\pi} e^{i\lambda(k-l)} f(\lambda) \, d\lambda$ for $1 \leq k, l \leq n$. Using the matrix identities $\det(AB^{-1}) = \det(I + B^{-1/2}(A - B)B^{-1/2})$ and $A^{-1} - B^{-1} = A^{-1}(A - B)B^{-1}$, we can write

$$\log \frac{p_f^{(n)}}{p_g^{(n)}} = -\frac{1}{2} \log \det(I + T_n(g)^{-1/2} T_n(f - g) T_n(g)^{-1/2})$$

$$- \frac{1}{2}(X^{(n)})^T T_n(f)^{-1} T_n(g - f) T_n(g)^{-1} X^{(n)}.$$

For a random vector $X$ with mean zero and covariance matrix $\Sigma$, we have $\mathrm{E}(X^T A X) = \mathrm{tr}(\Sigma A)$ and $\mathrm{var}(X^T A X) = \mathrm{tr}(\Sigma A \Sigma A) + \mathrm{tr}(\Sigma A \Sigma A^T)$. Hence,

$$P_f^{(n)}\left(\log \frac{p_f^{(n)}}{p_g^{(n)}}\right) = -\frac{1}{2}\log\det(I + T_n(g)^{-1/2}T_n(f-g)T_n(g)^{-1/2})$$

$$-\frac{1}{2}\mathrm{tr}(T_n(g-f)T_n(g)^{-1}),$$

$$4\,\mathrm{var}_{P_f^{(n)}}\left(\log \frac{p_f^{(n)}}{p_g^{(n)}}\right) = \mathrm{tr}(T_n(g-f)T_n(g)^{-1}T_n(g-f)T_n(g)^{-1})$$

$$+ \mathrm{tr}(T_n(g-f)T_n(g)^{-1}T_n(f)T_n(g)^{-1}T_n(g-f)T_n(f)^{-1}).$$

Define matrix norms by $\|A\|^2 = \sum_k \sum_l a_{k,l}^2 = \mathrm{tr}(AA^T)$ and $|A| = \sup\{\|Ax\| : \|x\| = 1\}$, where $\|x\|$ is the Euclidean norm. Then $\mathrm{tr}(A^2) \leq \|A\|^2$ and $\|AB\| \leq |A|\|B\|$. Furthermore, as a result of the inequalities $-\frac{1}{2}\mu^2 \leq \log(1+\mu) - \mu \leq 0$, for all $\mu \geq 0$, we have for any nonnegative definite matrix $A$ that $-\frac{1}{2}\mathrm{tr}(A^2) \leq \log\det(I + A) - \mathrm{tr}(A) \leq 0$. In view of the identities $x^T T_n(f) x = \int |\sum_k x_k e^{ik\lambda}|^2 f(\lambda)\, d\lambda$ and $x^T T_n(1) x = 2\pi \|x\|^2$, we also have that $|T_n(f)| \leq 2\pi \|f\|_\infty$ and $|T_n(f)^{-1}| \leq (2\pi)^{-1}\|1/f\|_\infty$. To see the validity of the second inequality, we use the fact that $\|A^{-1}\| \leq c^{-1}$ if $\|Ax\| \geq c\|x\|$ for all $x$. For $f \in \mathcal{F}$, $\|f\|_\infty < \infty$ and $\|1/f\|_\infty < \infty$. Furthermore,

$$\|T_n(f)\|^2 = \sum_{|h|<n}(n-|h|)\gamma_h^2(f) \leq 2\pi n \int_{-\pi}^{\pi} f^2(\lambda)\, d\lambda. \tag{6.2}$$

Using the preceding inequalities and the identity $\mathrm{tr}(AB) = \mathrm{tr}(BA)$, it is straightforward to obtain the desired bounds on the mean and variance of $\log(p_f^{(n)}/p_g^{(n)})$.

□

The preceding lemmas can be combined to obtain the following theorem, where the constants $\xi$ and $K$ are those introduced in Lemma 7.

THEOREM 7. *Let $P_f^{(n)}$ be the distribution of $(X_1, \ldots, X_n)$ for a stationary Gaussian time series $\{X_t : t = 0, \pm 1, \ldots\}$ with spectral density $f \in \mathcal{F}$. Assume that there exist constants $\Gamma$ and $M$ such that $\|\log f\|_\infty \leq \Gamma$ and $\sum_h |h|\gamma_h^2(f) \leq M$ for every $f \in \mathcal{F}$. Let $\varepsilon_n \geq 1/\sqrt{n}$ satisfy, for every $j \in \mathbb{N}$,*

$$\sup_{\varepsilon > \varepsilon_n} \log N(\xi\varepsilon/2, \{f \in \mathcal{F} : \|f - f_0\|_2 \leq \varepsilon\}, \|\cdot\|_\infty) \leq n\varepsilon_n^2,$$

$$\frac{\Pi(f : \|f - f_0\|_2 \leq j\varepsilon)}{\Pi(f : \|f - f_0\|_2 \leq \varepsilon)} \lesssim e^{Kn\varepsilon_n^2 j^2/8}.$$

*Then $P_{f_0}^{(n)}\Pi(f : \|f - f_0\|_2 \geq M_n \varepsilon_n | X_1, \ldots, X_n) \to 0$ for every $M_n \to \infty$.*

**7. Applications.** In this section, we present a number of examples of application of the general results obtained in the preceding sections. The examples concern combinations of a variety of models with various prior distributions.

7.1. *Finite sieves.* Consider the setting of independent, nonidentically distributed observations of Section 3. We construct sequences of priors, each supported on finitely many points such that the posterior distribution converges at a rate equal to the solution of an equation involving bracketing entropy numbers. Because bracketing entropy numbers are often close to metric entropy numbers, this construction exhibits priors for which the prior mass condition (2.5) is automatically satisfied. The construction is similar to that for the i.i.d. case given by Ghosal, Ghosh and van der Vaart [14], Section 3. However, in this case, some extra care is needed to appropriately define the bracketing numbers in the product space of densities. In the following, we consider a componentwise bracketing.

Consider a sequence of models $\mathcal{P}^{(n)} = \{P_\theta^{(n)} : \theta \in \Theta\}$ of $n$-fold product measures $P_\theta^{(n)}$, where each measure is given by a density $(x_1, \ldots, x_n) \mapsto \prod_{i=1}^n p_{\theta,i}(x_i)$ relative to a product-dominating measure $\bigotimes_{i=1}^n \mu_i$. For a given $n$ and $\varepsilon > 0$, we define the *componentwise Hellinger upper bracketing number* for $\Theta$ to be the smallest number $N$ such that there are integrable nonnegative functions $u_{j,i}$ for $j = 1, 2, \ldots, N$ and $i = 1, 2, \ldots, n$, with the property that for any $\theta \in \Theta$, there exists some $j$ such that $p_{\theta,i} \leq u_{j,i}$ for all $i = 1, 2, \ldots, n$ and $\sum_{i=1}^n h^2(p_{\theta,i}, u_{j,i})^2 \leq n\varepsilon^2$. We shall denote this by $N_]^{n\otimes}(\varepsilon, \Theta, d_n)$.

Given a sequence of sets $\Theta_n \uparrow \Theta$ and $\varepsilon_n \to 0$ such that $\log N_]^{n\otimes}(\varepsilon_n, \Theta_n, d_n) \leq n\varepsilon_n^2$, let $(u_{j,i} : j = 1, 2, \ldots, N, \ i = 1, 2, \ldots, n)$ be a componentwise Hellinger upper bracketing for $\Theta_n$ [where $N = N_]^{n\otimes}(\varepsilon_n, \Theta_n, d_n)$]. From this bracketing, we construct a prior distribution $\Pi_n$ on the collection of densities of product measures, by defining $\Pi_n$ to be the measure that assigns mass $N^{-1}$ to each of the joint densities $p_j^{(n)} = \otimes_{i=1}^n (u_{j,i}/\int u_{j,i}\, d\mu_i)$, $j = 1, 2, \ldots, N$. The collection $\mathcal{P}_n = \{p_j^{(n)} : j = 1, 2, \ldots, N\}$ forms a sieve for the models $\mathcal{P}^{(n)}$ and can be considered as the parameter space for a given $n$. Although it is possible for the spaces $\mathcal{P}_n$ to not be embedded in a fixed space, Theorem 4 still applies and implies the following result.

THEOREM 8. *Let $\Theta_n \uparrow \Theta$ and $\theta_0 \in \Theta$. Assume that $\log N_]^{n\otimes}(\varepsilon_n, \Theta_n, d_n) \leq n\varepsilon_n^2$ for some sequence $\varepsilon_n \to 0$ with $n\varepsilon_n^2 \to \infty$. Let $\Pi_n$ be the uniform measure on the renormalized collection of upper product brackets, as indicated previously. Then for all sufficiently large $M$,*

$$(7.1) \qquad P_{\theta_0}^{(n)} \Pi_n(p^{(n)} : d_n^2(p_{\theta_0}^{(n)}, p^{(n)}) \geq M\varepsilon_n^2 | X_1, X_2, \ldots, X_n) \to 0.$$

PROOF. As $\mathcal{P}_n$ consists of finitely many points, its covering number with respect to any metric is bounded by its cardinality. Thus, (3.2) holds and (3.3) holds trivially.

Let $\theta_0 \in \Theta_n$ for all $n > n_0$. For a given $n > n_0$, let $j_0$ be the index for which $p_{\theta_0,i} \leq u_{j_0,i}$ and $\sum_{i=1}^n h^2(p_{\theta_0,i}, u_{j_0,i}) \leq n\varepsilon_n^2$. If $p$ is a probability density, $u$ is an integrable function such that $u \geq p$ and $v = u/\int u$, then because $2ab \leq (a^2 + b^2)$, it easily follows that $h^2(p,v) \leq (\int u\, d\mu)^{-1/2} h^2(p,u)$.

For any two probability densities $p$ and $q$, we have (see, e.g., Lemma 8 of [17])

$$K(p,q) \lesssim h^2(p,q)\left(1 + \log\left\|\frac{p}{q}\right\|_\infty\right), \qquad V(p,q) \lesssim h^2(p,q)\left(1 + \log\left\|\frac{p}{q}\right\|_\infty\right)^2.$$

Together with the elementary inequalities $1 + \log x \leq 2\sqrt{x}$ and $(1 + \log x)^2 \leq (4x^{1/4})^2 = 16 x^{1/2}$ for all $x \geq 1$, the bounds imply that

$$K(p,q) \lesssim h^2(p,q)\left\|\frac{p}{q}\right\|_\infty^{1/2}, \qquad V(p,q) \lesssim h^2(p,q)\left\|\frac{p}{q}\right\|_\infty^{1/2}.$$

Because $(p_{\theta_0,i}/v_{j_0,i}) \leq \int u_{j_0,i}\, d\mu$, it follows that $n^{-1}\sum_{i=1}^n K(p_{\theta_0,i}, v_{j_0,i}) \lesssim \varepsilon_n^2$ and $n^{-1}\sum_{i=1}^n V(p_{\theta_0,i}, v_{j_0,i}) \lesssim \varepsilon_n^2$. Thus, $\prod_{i=1}^n v_{j_0,i}$ gets prior probability equal to $N^{-1} \geq e^{-n\varepsilon_n^2}$ and hence relation (3.4) also holds for a multiple of the present $\varepsilon_n$. Thus, the posterior converges at the rate $\varepsilon_n$ with respect to the metric $d_n$. $\square$

7.1.1. *Nonparametric Poisson regression.* Let $X_1, X_2, \ldots$ be independent Poisson-distributed random variables with parameters $\psi(z_1), \psi(z_2), \ldots$, where $\psi: \mathbb{R} \to (0, \infty)$ is an unknown increasing link function and $z_1, z_2, \ldots$ are one-dimensional covariates. We assume that $L \leq \psi \leq U$ for some constants $0 < L < U < \infty$.

If $l \leq \psi \leq u$, then for any $z$ and $x$, we have $e^{-\psi(z)}(\psi(z))^x/x! \leq e^{-l(z)}(u(z))^x/x!$. For a pair of link functions $l \leq u$, let $q_{l,u}(x,z) = e^{-l(z)}(u(z))^x/x!$ and put $f_{l,u}^{(n)}(x_1, x_2, \ldots, x_n) = \prod_{i=1}^n q_{l,u}(x_i, z_i)$. For any constants $L < \lambda_1, \lambda_2, \mu_1, \mu_2 < U$, we have

$$\sum_{x=0}^\infty \left(\left(e^{-\lambda_1}\frac{\mu_1^x}{x!}\right)^{1/2} - \left(e^{-\lambda_2}\frac{\mu_2^x}{x!}\right)^{1/2}\right)^2$$
$$= (e^{-(\lambda_1 + \mu_1)/2} - e^{-(\lambda_2 + \mu_2)/2})^2 + 2e^{-(\lambda_1 + \lambda_2)/2}(e^{(\mu_1 + \mu_2)/2} - e^{\sqrt{\mu_1\mu_2}})$$
$$\leq \left(\frac{1}{2} + \frac{1}{4}L^{-1}\right)e^{U-L}(|\lambda_1 - \lambda_2|^2 + |\mu_1 - \mu_2|^2).$$

Let $l_1 \leq u_1$ and $l_2 \leq u_2$ be two pairs of link functions taking their values in the interval $[L, U]$. Therefore, with $\mathbb{P}_n^z = n^{-1}\sum_{i=1}^n \delta_{z_i}$ being the empirical distributions of $z_1, z_2, \ldots, z_n$, we have that $d_n^2(f_{l_1,u_1}^{(n)}, f_{l_2,u_2}^{(n)}) \lesssim \int(|l_1 - l_2|^2 + |u_1 - u_2|^2) d\mathbb{P}_n^z$. Hence, an $\varepsilon$-bracketing of the link functions with respect to the $L_2(\mathbb{P}_n^z)$-metric yields a componentwise Hellinger upper bracketing whose size is a multiple of $\varepsilon$. Now the $\varepsilon$-bracketing entropy numbers of the above class are bounded by a multiple of $\varepsilon^{-1}$, relative to any $L_2$-metric (cf. Theorem 2.7.5 of [31]). Equating this with $n\varepsilon^2$, we obtain the rate $n^{-1/3}$ for posterior convergence, which is also the minimax rate, relative to $d_n$.

In this example, the normalized upper brackets for the densities are also Poisson mass functions corresponding to the link functions equal to the upper brackets. Hence, the prior can be viewed as charging the space of link functions and the distance $d_n$ can also be induced on this space. This makes interpretations of the prior and the posterior, as well as the posterior convergence rate, more transparent. Further, as the space of link functions is a fixed space, proceeding as in Theorem 3.1 of [14], a fixed prior not depending on $n$ may be constructed such that the posterior converges at the same $n^{-1/3}$ rate.

7.2. *Linear regression with unknown error distribution.* Let $X_1, \ldots, X_n$ be independent regression response variables satisfying $X_i = \alpha + \beta z_i + \varepsilon_i$, $i = 1, 2, \ldots, n$, where the $z_i$'s are nonrandom one-dimensional covariates lying in $[-L, L]$ for some $L$ and the errors $\varepsilon_i$ are i.i.d. with density $f$ following some prior $\Pi$. Amewou-Atisso et al. [1] studied posterior consistency under this setup. Here, we refine the result to a posterior convergence rate. Assume that $|f'(x)| \leq C$ for all $x$ and all $f$ in the support of $\Pi$. The priors for $\alpha$ and $\beta$ are assumed to be compactly supported with positive densities in the interiors of their supports and all the parameters are assumed to be a priori independent. Let the true value of $(f, \alpha, \beta)$ be $(f_0, \alpha_0, \beta_0)$, an interior point in the support of the prior.

Let $H(\varepsilon)$ be a bound for the Hellinger $\varepsilon$-entropy of the support of $\Pi$ and suppose that $f_0(x)/f(x) \leq M(x)$ for all $x$, where $\int M^\delta f_0 < \infty$, $\delta > 0$. Then by Theorem 5 of [33], it follows that $\max\{K(f_0, f), V(f_0, f)\} \lesssim h^2(f_0, f) \times \log^2(1/h(f_0, f))$. Let $a(\varepsilon) = -\log \Pi(h(f_0, f) \leq \varepsilon)$. The posterior convergence rate for density estimation is then $\varepsilon_n$, given by

(7.2) $$\max\{H(\varepsilon_n), a(\varepsilon_n/(\log \varepsilon_n^{-1}))\} \leq n\varepsilon_n^2.$$

The following theorem shows that Euclidean parameters do not affect the rate.

THEOREM 9. *Under the above setup, if $f_0(x - \alpha_0 - \beta_0 z)/f(x - \alpha - \beta z) \leq M(x)$ for all $x, z, \alpha, \beta$, then the joint posterior of $(\alpha, \beta, f)$ concentrates around $(\alpha_0, \beta_0, f_0)$ at the rate $\varepsilon_n$ defined by (7.2), with respect to $d_n$.*

PROOF. We have, by $(a+b)^2 \leq 2(a^2+b^2)$, that $h^2(f_1(\cdot-\alpha_1-\beta_1 z), f_2(\cdot-\alpha_2-\beta_2 z)) \leq 2h^2(f_1, f_2) + 4C^2|\alpha_1-\alpha_2|^2 + 4C^2L^2|\beta_1-\beta_2|^2$, which leads to

$$d_n^2(P_{f_1,\alpha_1,\beta_1}^{(n)}, P_{f_2,\alpha_2,\beta_2}^{(n)}) \lesssim h^2(f_1, f_2) + |\alpha_1-\alpha_2|^2 + |\beta_1-\beta_2|^2$$

and hence the $d_n$-entropy of the parameter space is bounded by a multiple of $H(\varepsilon) + \log\frac{1}{\varepsilon} \lesssim H(\varepsilon)$.

To lower bound the prior probability of $\bar{B}_n((f_0, \alpha_0, \beta_0), \varepsilon; 2)$ defined by (3.5), by Theorem 5 of [33] with $h = h(f_0(\cdot - \alpha_0 - \beta_0 z), f(\cdot - \alpha - \beta z))$, we have that $K(f_0(\cdot - \alpha_0 - \beta_0 z), f(\cdot - \alpha - \beta z)) \lesssim h^2 \log\frac{1}{h}$ and $V(f_0(\cdot - \alpha_0 - \beta_0 z), f(\cdot - \alpha - \beta z)) \lesssim h^2 \log^2\frac{1}{h}$. Thus, a multiple of $\varepsilon^{-2} e^{-ca(\varepsilon/\log \varepsilon^{-1})}$ lower bounds the prior probability of (3.5) and the first factor can be absorbed into the second, where $c$ is a suitable positive constant. Thus, Theorem 4 implies that the posterior convergence rate with respect to $d_n$ is $\varepsilon_n$. □

More concretely, if the prior is a Dirichlet mixture of normals (or its symmetrization) with the scale parameter lying between two positive numbers and the base measure having compact support, and if the true error density is also a normal mixture of this type, then by Ghosal and van der Vaart [16], it follows that the convergence rate is $(\log n)/\sqrt{n}$. The assumption of compact support of the base measure can be relaxed by using sieves. Compactness of the support of the prior for $\alpha$ and $\beta$ may be relaxed by using sieves $|\alpha| \leq c \log n$ if these priors have sub-Gaussian tails. Also, it is straightforward to extend the result to a multidimensional regressor. For more general error densities, one has to allow arbitrarily small scale parameters and apply the results of Ghosal and van der Vaart [17] to obtain a slower rate.

Often, only the Euclidean part is of interest and an $n^{-1/2}$ rate of convergence is generally obtained in the classical context. The posterior of the Euclidean part is also expected to converge at an $n^{-1/2}$ rate and the Bernstein–von Mises theorem may hold; see [26] for some results. However, as we consider $(f, \alpha, \beta)$ together and obtain global convergence rates, it seems unlikely that our methods will yield these improved convergence rates for the Euclidean portion of the parameter.

7.3. *Whittle estimation of the spectral density.* Let $\{X_t : t \in \mathbb{Z}\}$ be a second order stationary time series with mean zero and autocovariance function $\gamma_r = \mathrm{E}(X_t X_{t+r})$. The spectral density of the process is defined (under the assumption that $\sum_r |\gamma_r| < \infty$) by $f(\lambda) = \frac{1}{2\pi} \sum_{r=-\infty}^{\infty} \gamma_r e^{-ir\pi\lambda}$, $\lambda \in [0, 1]$; here, we have changed the original domain $[-\pi, \pi]$ of spectral density to $[0, 1]$ by using symmetry and then rescaling. Let $I_n(\lambda) = (2\pi n)^{-1} |\sum_{t=1}^{n} X_t e^{-it\pi\lambda}|^2$, $\lambda \in [0, 1]$, denote the periodogram. Because the likelihood is complicated, Whittle [32] proposed as an approximate likelihood that of a sample $U_1, \ldots, U_\nu$ of independent exponential variables with means $f(2j/n)$, $j = 1,\ldots,\nu$, evaluated with $U_j = I_n(2j/n)$, where $\nu = \lfloor n/2 \rfloor$. The Whittle likelihood is motivated by the fact that if $\lambda_{n,i} \to \lambda_i$, $i = 1,\ldots,m$, then under reasonable conditions such as mixing conditions, $(I_n(\lambda_{n,1}),\ldots,I_n(\lambda_{n,m}))$ converges weakly to a vector of independent exponential variables with mean vector $(f(\lambda_1),\ldots,f(\lambda_m))$; see, for instance, Theorem 10.3.2 of Brockwell and Davis [6]. Dahlhaus [10] applied the technique of Whittle likelihood to estimating the spectral density by the minimum contrast method. A consistent Bayesian nonparametric method has been proposed by Choudhuri, Ghosal and Roy [7]. Below, we indicate how to obtain a rate of convergence using Theorem 4.

As in the proof of consistency, we use the contiguity result of Choudhuri, Ghosal and Roy [8], which shows that for a Gaussian time series, the sequence of laws of $(I_n(2/n),\ldots,I_n(2\nu/n))$ and the sequence of approximating exponential distributions of $(U_1,\ldots,U_\nu)$ are contiguous. Thus, a rate of convergence of the posterior distribution under the actual distribution follows from a rate of convergence under the assumption that $U_1,\ldots,U_\nu$ are exactly independent and exponentially distributed with means $f(2/n),\ldots,f(2\nu/n)$, to which Theorem 4 can be applied.

Let $\bar{d}_n^2(f_1, f_2) = \nu^{-1} \sum_{i=1}^{\nu} (f_1(2i/n) - f_2(2i/n))^2$. If $f_1$ and $f_2$ are spectral densities with $m \le f_1, f_2 \le M$ pointwise, then it follows that

$$(7.3) \quad \frac{1}{4M^2} \bar{d}_n^2(f_1, f_2) \le d_n^2(f_1, f_2) \le \frac{1}{4m^2} \bar{d}_n^2(f_1, f_2) \le \frac{1}{4m^2} \|f_1 - f_2\|_\infty^2,$$

where $d_n$ is given by (3.1) and $\|\cdot\|_\infty$ is the uniform distance. If the spectral densities are Lipschitz continuous, then a rate for the discretized $L_2$-distance $d_n$ will imply a rate for the ordinary $L_2$-distance $\|\cdot\|_2$ by the relation $\|f_1 - f_2\|_2 \lesssim \bar{d}_n(f_1, f_2) + (L + M)/n$, where $L$ and $M$ are the Lipschitz constant and uniform bound, respectively. To see this, note that $\sqrt{n/\nu}\bar{d}_n(f, 0) = \|f_n\|_2$, where $f_n = \sum_{j=1}^{\nu} f(2j/n) 1_{((2j-2)/n, 2j/n]}$ and hence

$$|\|f\|_2 - \sqrt{n/\nu}\bar{d}_n(f, 0)| \lesssim \|f - f_n\|_2 \lesssim \frac{\|f\|_{\mathrm{Lip}}}{n} + \left(1 - \frac{2\nu}{n}\right)\|f\|_\infty.$$

It follows that for the verification of (3.2), we may always replace $d_n$ by $\bar{d}_n$ and if the spectral densities are restricted to Lipschitz functions with Lipschitz constant $L_n$ and where $\varepsilon_n \gg L_n/n$, then we may also replace $d_n$ by the $L_2$-norm $\|\cdot\|_2$.

Now, by easy calculations, for all spectral densities $f, f_0$ taking values in $[m, M]$, we have that $\nu^{-1} \sum_{i=1}^{\nu} K(P_{f_0,i}, P_{f,i}) \lesssim \bar{d}_n^2(f_0, f) \lesssim \|f - f_0\|_\infty^2$ and $\nu^{-1} \sum_{i=1}^{\nu} V_{2,0}(P_{f_0,i}, P_{f,i}) \lesssim \bar{d}_n^2(f_0, f) \lesssim \|f - f_0\|_\infty^2$, hence it suffices to estimate the prior probability of sets of the form $\{f : \|f - f_0\|_\infty \le \varepsilon\}$. Alternatively, if the spectral densities under consideration are Lipschitz, then we may estimate the prior mass of an $L_2$-ball around $f_0$.

As a concrete prior, we consider the prior used by Choudhuri, Ghosal and Roy [7], namely $f = \tau q$, where $\tau = \text{var}(X_t)$ has a nonsingular prior density and $q$, a probability density on $[0,1]$, is given the Dirichlet–Bernstein prior of Petrone [24]. We then restrict the prior to the set $\mathbb{K} = \{f : m < f < M\}$. The order of the Bernstein polynomial, $k$, has prior mass function $\rho$, which is assumed to satisfy $e^{-\beta_1 k \log k} \lesssim \rho(k) \lesssim e^{-\beta_2 k}$. Let $\Pi$ denote the resulting prior.

Clearly, as $f_0 \in \mathbb{K}$, restricting the prior to $\mathbb{K}$ can only increase the prior probability of $\{f : \|f - f_0\|_\infty < \varepsilon\}$. Therefore, following Ghosal [12], $\Pi(\|f - f_0\|_\infty < \varepsilon) \gtrsim e^{-c\varepsilon^{-1} \log \varepsilon^{-1}}$. Hence, $\varepsilon_n$ of the order $n^{-1/3}(\log n)^{1/3}$ satisfies (3.4).

Consider a sieve $\mathcal{F}_n$ for the parameter space $\mathbb{K}$, which consists solely of Bernstein polynomials of order $k_n$ or less. All of these functions have Lipschitz constant at most $k_n^2$ and are uniformly bounded away from zero and infinity by construction. The $\varepsilon$-entropy of $\mathcal{F}_n$ relative to $\bar{d}_n$ can be bounded above by that of the simplex, which is further bounded above by $k \log k + k \log \varepsilon^{-1}$. Hence, by choosing $k_n$ of the order $n^{1/3}(\log n)^{2/3}$, the convergence rate at $f_0$ on $\mathcal{F}_n$ with respect to $d_n$ is given by $\max(n^{-1/2} k_n^{1/2}(\log n)^{1/2}, n^{-1/3}(\log n)^{1/3}, k_n^2/n) = n^{-1/3}(\log n)^{4/3}$. Now, $\Pi(\mathcal{F}_n^c) = \rho(k > k_n) \lesssim e^{-\beta_2 k_n} = e^{-\beta n^{1/3}(\log n)^{2/3}} = e^{-\beta n(n^{-1/3}(\log n)^{1/3})^2}$. Thus, the posterior probability of $\mathcal{F}_n^c$ goes to zero by Lemma 1 and hence the convergence rate on $\mathbb{K}$ is also $n^{-1/3}(\log n)^{1/3}$. The minimax rate $n^{-2/5}$ may be obtained, for instance, by using splines, which have better approximation properties.

7.4. *Nonlinear autoregression.* Consider the nonlinear autoregressive model in which we observe the elements $X_1, \ldots, X_n$ of a stationary time series $\{X_t : t \in \mathbb{Z}\}$ satisfying

$$(7.4) \qquad X_i = f(X_{i-1}) + \varepsilon_i, \qquad i = 1, 2, \ldots, n,$$

where $f$ is an unknown function and $\varepsilon_1, \varepsilon_2, \ldots, \varepsilon_n$ are i.i.d. $N(0, \sigma^2)$. For simplicity, we assume that $\sigma = 1$. Then $X_n$ is a Markov chain with transition density $p_f(y|x) = \phi(y - f(x))$, where $\phi(x) = (2\pi)^{-1/2} e^{-x^2/2}$. Assume that $f \in \mathcal{F}$, a class of functions such that $|f(x)| \leq M$ and $|f(x) - f(y)| \leq L|x - y|$ for all $x, y$ and $f \in \mathcal{F}$.

Set $r(y) = \frac{1}{2}(\phi(y - M) + \phi(y + M))$. Then $r(y) \lesssim p_f(y|x) \lesssim r(y)$ for all $x, y \in \mathbb{R}$ and $f \in \mathcal{F}$. Further, $\sup\{\int |p(y|x_1) - p(y|x_2)| \, dy : x_1, x_2 \in \mathbb{R}\} < 2$. Hence, the chain is $\alpha$-mixing with exponentially decaying mixing coefficients and has a unique stationary distribution $Q_f$ whose density $q_f$ satisfies $r \lesssim q_f \lesssim r$. Let $\|f\|_s = (\int |f|^s \, dr)^{1/s}$.

Because $h^2(N(\mu_1, 1), N(\mu_2, 1)) = 2[1 - \exp(-|\mu_1 - \mu_2|^2/8)]$, it easily follows for $f_1, f_2 \in \mathcal{F}$, $d$ defined in (4.2) and $d\nu = r \, d\lambda$ that $\|f_1 - f_2\|_2 \lesssim$

$d(f_1, f_2) \lesssim \|f_1 - f_2\|_2$. Thus, we may verify (4.5) relative to the $L_2(r)$-metric. It can also be computed that

$$P_{f_0} \log \frac{p_{f_0}(X_2|X_1)}{p_f(X_2|X_1)} = \frac{1}{2} \int (f_0 - f)^2 q_{f_0} \, d\lambda \lesssim \|f - f_0\|_2^2,$$

$$P_{f_0} \left| \log \frac{p_{f_0}(X_2|X_1)}{p_f(X_2|X_1)} \right|^s \lesssim \int |f_0 - f|^s q_{f_0} \, d\lambda \lesssim \|f - f_0\|_s^s.$$

Thus, $B^*(f_0, \varepsilon; s) \supset \{f : \|f - f_0\|_s \leq c\varepsilon\}$ for some constant $c > 0$, where $B^*(f_0, \varepsilon; s)$ is as in Theorem 5. Thus, it suffices to verify (4.7) with $s > 2$.

7.4.1. *Random histograms.* As a prior on the regression functions $f$, consider a random histogram as follows. For a given number $K \in \mathbb{N}$, partition a given compact interval in $\mathbb{R}$ into $K$ intervals $I_1, \ldots, I_K$ and let $I_0 = \mathbb{R} \setminus \bigcup_k I_k$. Let the prior $\Pi_n$ on $f$ be induced by the map $\alpha \mapsto f_\alpha$ given by $f_\alpha = \sum_{k=1}^K \alpha_k 1_{I_k}$, where the coordinates $\alpha_1, \ldots, \alpha_K$ of $\alpha \in \mathbb{R}^K$ are chosen to be i.i.d. random variables with the uniform distribution on the interval $[-M, M]$ and where $K = K_n$ is to be chosen later. Let $r(I_k) = \int_{I_k} r \, d\lambda$.

The support of $\Pi_n$ consists of all functions with values in $[-M, M]$ that are piecewise constant on each interval $I_k$ for $k = 1, \ldots, K$ and which vanish on $I_0$. For any pair $f_\alpha$ and $f_\beta$ of such functions, we have, for any $s \in [2, \infty]$, $\|f_\alpha - f_\beta\|_s = \|\alpha - \beta\|_s$, where $\|\alpha\|_s$ is the $r$-weighted $\ell_s$-norm of $\alpha = (\alpha_1, \ldots, \alpha_K) \in \mathbb{R}^K$ given by $\|\alpha\|_s^s = \sum_k |\alpha_k|^s r(I_k)$. The dual use of $\|\cdot\|_s$ should not lead to any confusion as it will be clear from the context whether $\|\cdot\|_s$ is a norm on functions or on vectors. The $L_2(r)$-projection of $f_0$ onto this support is the function $f_{\alpha_0}$ for $\alpha_{0,k} = \int_{I_k} f_0 r \, d\lambda / r(I_k)$, whence, by Pythagoras' theorem, $\|f_\alpha - f_0\|_2^2 = \|f_\alpha - f_{\alpha_0}\|_2^2 + \|f_{\alpha_0} - f_0\|_2^2$ for any $\alpha \in [-M, M]^K$. In particular, $\|f_\alpha - f_0\|_2 \geq c\|\alpha - \alpha_0\|_2$ for some constant $c$ and hence, with $\mathcal{F}_n$ denoting the support of $\Pi_n$,

$$N(\varepsilon, \{f \in \mathcal{F}_n : \|f - f_0\|_2 \leq 16\varepsilon\}, \|\cdot\|_2)$$
$$\leq N(\varepsilon, \{\alpha \in \mathbb{R}^K : \|\alpha - \alpha_0\|_2 \leq 16c\varepsilon\}, \|\cdot\|_2) \leq (80c)^K,$$

as in Lemma 4.1 of [25]. Thus, (4.5) holds if $n\varepsilon_n^2 \gtrsim K$.

To verify (4.7), note that for $\lambda = (\lambda(I_1), \ldots, \lambda(I_K))$,

$$\|f_{\alpha_0} - f_0\|_s^s = \int_{I_0} |f_0|^s \, d\lambda + \sum_k \int_{I_k} |\alpha_{0,k} - f_0|^s r \, d\lambda \leq M^s r(I_0) + L^s \|\lambda\|_s^s.$$

Hence, as $f_0 \in \mathcal{F}$, for every $\alpha \in [-M, M]^K$,

$$\|f_\alpha - f_0\|_s \lesssim \|\alpha - \alpha_0\|_s + r(I_0)^{1/s} + \|\lambda\|_s \leq \|\alpha - \alpha_0\|_\infty + r(I_0)^{1/s} + \|\lambda\|_s,$$

where $\|\cdot\|_\infty$ is the ordinary maximum norm on $\mathbb{R}^K$. For $r(I_0)^{1/s} + \|\lambda\|_s \leq \varepsilon/2$, we have that $\{f : \|f - f_0\|_s \leq \varepsilon\} \supset \{f_\alpha : \|\alpha - \alpha_0\|_\infty \leq \varepsilon/2\}$. Using $\|\alpha - \alpha_0\|_2 \leq c\|f_\alpha - f_0\|_2$, for any $\varepsilon > 0$ such that $r(I_0)^{1/s} + \|\lambda\|_s \leq \varepsilon/2$, we have

$$\frac{\Pi_n(f : \|f - f_0\|_2 \leq j\varepsilon)}{\Pi_n(f : \|f - f_0\|_s \leq \varepsilon)} \leq \frac{\Pi_n(\alpha : \|\alpha - \alpha_0\|_2 \leq j\varepsilon)}{\Pi_n(\alpha : \|\alpha - \alpha_0\|_\infty \leq \varepsilon c/2)}.$$

We show that the right-hand side is bounded by $e^{Cn\varepsilon^2/8}$ for some $C$.

For $\bigcup_k I_k$, a regular partition of an interval $[-A, A]$, we have that $\|\lambda\|_s = 2A/K$ and since $r(I_k) \geq \lambda(I_k) \inf_{x \in I_k} r(x)$ for every $k \geq 1$, the norm $\|\cdot\|_2$ is bounded below by $\sqrt{2A\phi(A)/K} \gtrsim \sqrt{\phi(A)/K}$ times a multiple of the Euclidean norm. In this case, the preceding display is bounded above by

$$\frac{(Cj\varepsilon\sqrt{K/\phi(A)}/(2M))^K \text{vol}_K}{(\varepsilon c/(4M))^K} \sim \left(\frac{j\sqrt{2\pi e}}{\sqrt{\phi(A)}}\right)^K \frac{1}{\sqrt{\pi K}},$$

by Stirling's approximation, where $\text{vol}_K$ is the volume of the $K$-dimensional Euclidean unit ball. The probability $r(I_0)$ is bounded above by $1 - 2\Phi(A) \lesssim \phi(A)$. Hence, (4.7) will hold if $K\log(1/\phi(A)) \lesssim n\varepsilon_n^2$, $\phi(A) \lesssim \varepsilon_n^s$ and $A/K \lesssim \varepsilon_n$. All requirements are met for $\varepsilon_n$ equal to a multiple of $n^{-1/3}(\log n)^{1/2}$ [with $K \sim \sqrt{\log(1/\varepsilon_n)}\varepsilon_n^{-1}$ and $A \sim \sqrt{\log(1/\varepsilon_n)}$]. This is only marginally weaker than the minimax rate, which is $n^{-1/3}$ for this problem, provided the autoregression functions are assumed to be only Lipschitz continuous.

The logarithmic factor in the convergence rate appears to be a consequence of the fact that the regression functions are defined on the full real line. The present prior is a special case of a spline-based prior (see, e.g., Section 7.7). If $f$ has smoothness beyond Lipschitz continuity, then the use of higher order splines should yield a faster convergence rate.

7.5. *Finite-dimensional i.n.i.d. models.* Theorem 4 is also applicable to finite-dimensional models and yields the usual convergence rate as shown below. The result may be compared with Theorem I.10.2 of [19] and Proposition 1 of [13].

THEOREM 10. *Let $X_1, \ldots, X_n$ be i.n.i.d. observations following densities $p_{\theta,i}$, where $\Theta \subset \mathbb{R}^d$. Let $\theta_0$ be an interior point of $\Theta$. Assume that there exist constants $\alpha > 0$ and $0 \leq c_i \leq C_i < \infty$ with, for every $\theta, \theta_1, \theta_2 \in \Theta$,*

$$(7.5) \qquad c = \liminf_{n \to \infty} \frac{1}{n} \sum_{i=1}^n c_i > 0, \qquad C = \limsup_{n \to \infty} \frac{1}{n} \sum_{i=1}^n C_i < \infty$$

*such that $P_{\theta_0,i}(\log \frac{p_{\theta_0,i}}{p_{\theta,i}}) \leq C_i\|\theta - \theta_0\|^{2\alpha}$, $P_{\theta_0,i}(\log \frac{p_{\theta_0,i}}{p_{\theta,i}})^2 \leq C_i\|\theta - \theta_0\|^{2\alpha}$ and*

$$(7.6) \qquad c_i\|\theta_1 - \theta_2\|^{2\alpha} \leq h^2(p_{\theta_1,i}, p_{\theta_2,i}) \leq C_i\|\theta_1 - \theta_2\|^{2\alpha}.$$

POSTERIOR CONVERGENCE RATES 21*Assume that the prior measure $\Pi$ possesses a density $\pi$ which is bounded away from zero in a neighborhood of $\theta_0$ and bounded above on the entire parameter space. Then the posterior converges at the rate $n^{-1/(2\alpha)}$ with respect to the Euclidean metric.*

For regular families, the above displays are satisfied for $\alpha = 1$ and the usual $n^{-1/2}$ rate is obtained; see [19], Chapter III. Nonregular cases, for instance, when the densities have discontinuities depending on the parameter [such as the uniform distribution on $(0, \theta)$], have $\alpha < 1$ and faster rates are obtained; see [19], Chapters V and VI and [13].

PROOF OF THEOREM 10. By the assumptions (7.5) and (7.6), it suffices to show that the posterior convergence rate with respect to $d_n$ defined by (3.1) is $n^{-1/2}$. Now, by Pollard ([25], Lemma 4.1),

$$N(\varepsilon/18, \{\theta \in \Theta : d_n(\theta, \theta_0) < \varepsilon\}, d_n)$$

(7.7)
$$\leq N((\varepsilon^2/(36C))^{1/(2\alpha)}, \{\theta \in \Theta : \|\theta - \theta_0\| < (2\varepsilon^2/c)^{1/(2\alpha)}\}, \|\cdot\|)$$

$$\leq 6^d \left(\frac{72C}{c}\right)^{d/(2\alpha)},$$

which verifies (3.2). For (3.4), note that

$$\frac{\Pi(\theta : d_n(p_\theta, p_{\theta_0}) \leq j\varepsilon)}{\Pi(\theta : n^{-1}\sum_{i=1}^n K_i(\theta_0, \theta) \leq \varepsilon^2, n^{-1}\sum_{i=1}^n V_{2;i}(\theta_0, \theta) \leq \varepsilon^2)}$$

$$\leq \frac{\Pi(\theta : \|\theta - \theta_0\| \leq (2j^2\varepsilon^2/c)^{1/(2\alpha)})}{\Pi(\theta : \|\theta - \theta_0\| \leq (\varepsilon^2/(2C))^{1/(2\alpha)})} \leq A j^{d/\alpha}$$

for sufficiently small $\varepsilon > 0$, where $A$ is a constant depending on $d$, $c$, $C$ and the upper and lower bounds on the prior density. The conclusion follows for $\varepsilon_n = M/\sqrt{n}$, where $M$ is a large constant. □

The condition that the Hellinger distance is bounded below by a power of the Euclidean distance excludes the possibility of unbounded parameter spaces. This defect may be rectified by applying Theorem 3 to derive the rate. If there is a uniformly exponentially consistent test for $\theta = \theta_0$ against the complement of a bounded set, then the result holds even if $\Theta$ is not bounded. Often, such tests exist by virtue of bounds on log affinity, as in the case of normal distributions, or by large deviation type inequalities; see [20] and [14], Section 7. Further, if the prior density is not bounded above, but has a polynomial or subexponential majorant, then the rate calculation also remains valid.

7.6. *White noise with conjugate priors.* In this section, we consider the white noise model of Section 5 with a conjugate Gaussian prior. This allows us to complement and rederive results of Zhao [34] and Shen and Wasserman [27] in our framework. Thus, we observe an infinite sequence $X_1, X_2, \ldots$ of independent random variables, where $X_i$ is normally distributed with mean $\theta_i$ and variance $n^{-1}$.

We consider the prior $\Pi_n$ on the parameter $\theta = (\theta_1, \theta_2, \ldots)$ that can be structurally described by saying that $\theta_1, \ldots, \theta_k$ are independent with $\theta_i$ normally distributed with mean zero and variance $\sigma_{i,k}^2$ and that $\theta_{k+1}, \theta_{k+2}, \ldots$ are set equal to zero. Here, we choose the cutoff $k$ dependent on $n$ and equal to $k = \lfloor n^{1/(2\alpha+1)} \rfloor$ for some $\alpha > 0$. Zhao [34] and Shen and Wasserman [27] consider the case where $\sigma_{i,k}^2 = i^{-(2\alpha+1)}$ for $i = 1, \ldots, k$ and show that the convergence rate is $\varepsilon_n = n^{-\alpha/(2\alpha+1)}$ if the true parameter $\theta_0$ is "$\alpha$-regular" in the sense that $\sum_{i=1}^\infty \theta_{0,i}^2 i^{2\alpha} < \infty$. We shall obtain the same result for any triangular array of variances such that

(7.8) $$\min\{\sigma_{i,k}^2 i^{2\alpha} : 1 \leq i \leq k\} \sim k^{-1}.$$

For instance, for each $k$, the coefficients $\theta_1, \ldots, \theta_k$ could be chosen i.i.d. normal with mean zero and variance $k^{-1}$ or could follow the model of the authors mentioned previously.

THEOREM 11. *If $k \sim n^{1/(2\alpha+1)}$ and (7.8) holds, then the posterior converges at the rate $\varepsilon_n = n^{-\alpha/(2\alpha+1)}$ for any $\theta_0$ such that $\sum_{i=1}^\infty \theta_{0,i}^2 i^{2\alpha} < \infty$.*

PROOF. The support $\Theta_n$ of the prior is the set of all $\theta \in \ell_2$ with $\theta_i = 0$ for $i > k$ and can be identified with $\mathbb{R}^k$. Moreover, the $\ell_2$-norm $\|\cdot\|$ on the support can be identified with the Euclidean norm $\|\cdot\|_k$ on $\mathbb{R}^k$. Let $B_k(x, \varepsilon)$ denote the $k$-dimensional Euclidean ball of radius $\varepsilon$ and center $x \in \mathbb{R}^d$. For any true parameter $\theta_0 \in \ell_2$, we have $\|\theta - \theta_0\| \geq \|\mathbf{P}\theta - \mathbf{P}\theta_0\|_k$, where $\mathbf{P}$ is the projection on $\Theta_n$, and hence

$$N(\varepsilon/8, \{\theta \in \Theta_n : \|\theta - \theta_0\| \leq \varepsilon\}, \|\cdot\|) \leq N(\varepsilon/8, B_k(\mathbf{P}\theta_0, \varepsilon), \|\cdot\|_k) \leq (40)^k.$$

It follows that (5.1) is satisfied for $n\varepsilon_n^2 \gtrsim k$, that is, in view of our choice of $k$, $\varepsilon_n \gtrsim n^{-\alpha/(2\alpha+1)}$.

By Pythagoras' theorem, we have that $\|\theta - \theta_0\|^2 = \|\mathbf{P}\theta - \mathbf{P}\theta_0\|^2 + \sum_{i>k} \theta_{0,i}^2$ for any $\theta$ in the support of $\Pi_n$. Hence, for $\sum_{i>k} \theta_{0,i}^2 \leq \varepsilon_n^2/2$, we have that

$$\Pi_n(\theta \in \Theta_n : \|\theta - \theta_0\| \leq \varepsilon_n) \geq \Pi_n(\theta \in \mathbb{R}^k : \|\theta - \mathbf{P}\theta_0\|_k \leq \varepsilon_n/2).$$

By the definition of the prior, the right-hand side involves a quadratic form in Gaussian variables. For $\Sigma$ the $k \times k$ diagonal matrix with elements $\sigma_{i,k}^2$, the quotient on the left-hand side of (5.2) can be bounded as

$$\frac{\Pi_n(\theta \in \Theta_n : \|\theta - \theta_0\| \leq j\varepsilon_n)}{\Pi_n(\theta \in \Theta_n : \|\theta - \theta_0\| \leq \varepsilon_n)} \leq \frac{N_k(-\mathbf{P}\theta_0, \Sigma)(B(0, j\varepsilon_n))}{N_k(-\mathbf{P}\theta_0, \Sigma)(B(0, \varepsilon_n/2))}.$$

The probability in the numerator increases if we center the normal distribution at 0 rather than at $-\mathbf{P}\theta_0$, by Anderson's lemma. Furthermore, for any $\mu \in \mathbb{R}^k$,

$$\frac{dN_k(\mu, \Sigma)}{dN_k(0, \Sigma/2)}(\theta) = \frac{e^{-\sum_{i=1}^k (\theta_i - \mu_i)^2/(2\sigma_{i,k}^2)}}{\sqrt{2}^k e^{-\sum_{i=1}^k \theta_i^2/\sigma_{i,k}^2}} \geq 2^{-k/2} e^{-\sum_{i=1}^k \mu_i^2/\sigma_{i,k}^2}.$$

Therefore, we may recenter the denominator at 0 at the cost of adding the factor on the right (with $\mu = \theta_0$) and dividing the covariance matrix by 2. We obtain that the left-hand side of (5.2) is bounded above by

$$2^{k/2} e^{\sum_{i=1}^k \theta_{0,i}^2/\sigma_{i,k}^2} \frac{N_k(0, \Sigma)(B(0, j\varepsilon_n))}{N_k(0, \Sigma/2)(B(0, \varepsilon_n/2))}$$

$$\leq 2^{k/2} e^{\sum_{i=1}^k \theta_{0,i}^2/\sigma_{i,k}^2} \left(\frac{\bar{\sigma}_k}{\underline{\sigma}_k}\right)^k \frac{N_k(0, \bar{\sigma}_k^2 I)(B(0, j\varepsilon_n))}{N_k(0, \underline{\sigma}_k^2 I/2)(B(0, \varepsilon_n/2))},$$

where $\bar{\sigma}_k$ and $\underline{\sigma}_k$ denote the maximum and the minimum of $\sigma_{i,k}$ for $i = 1, 2, \ldots, k$. The probabilities on the right-hand side are left tail probabilities of chi-square distributions with $k$ degrees of freedom, and can be expressed as integrals. The preceding display is bounded above by

$$2^{k/2} e^{\sum_{i=1}^k \theta_{0,i}^2/\sigma_{i,k}^2} \left(\frac{\bar{\sigma}_k}{\underline{\sigma}_k}\right)^k \frac{\int_0^{j^2 \varepsilon_n^2/\bar{\sigma}_k^2} x^{k/2-1} e^{-x/2} \, dx}{\int_0^{\varepsilon_n^2/(2\underline{\sigma}_k^2)} x^{k/2-1} e^{-x/2} \, dx}.$$

The exponential in the integral in the numerator is bounded above by 1 and hence this integral is bounded above by $j^k \varepsilon_n^k/(k\bar{\sigma}_k^k)$. We now consider two separate cases. If $\varepsilon_n^2/\underline{\sigma}_k^2$ remains bounded, then we can also bound the exponential in the integral in the denominator below by a constant and have that the preceding display is bounded above by a multiple of $4^k j^k \exp(\sum_{i=1}^k \theta_{0,i}^2/\sigma_{i,k}^2)$. If $\varepsilon_n^2/\underline{\sigma}_k^2 \to \infty$, then we bound the integral in the denominator below by $(\eta/2)^{k/2-1} \int_{\eta/2}^\eta e^{-x/2} \, dx$ for $\eta = \varepsilon_n^2/(2\underline{\sigma}_k^2)$. This leads to the upper bound being a multiple of $8^k j^k \exp(\sum_{i=1}^k \theta_{0,i}^2 \sigma_{i,k}^{-2}) \varepsilon_n^2 \underline{\sigma}_k^{-2} \exp(\varepsilon_n^2 \underline{\sigma}_k^{-2}/8)$. By the assumption (7.8), we have that $\underline{\sigma}_k^2 \gtrsim k^{-(2\alpha+1)} \sim n^{-1}$. We also have that $k \sim n\varepsilon_n^2$. It follows that $\varepsilon_n^2/\underline{\sigma}_k^2 \lesssim n\varepsilon_n^2$ and that $\underline{\sigma}_k^{-2}$ is bounded by a polynomial in $k$. We conclude that with our choice of $k \sim n^{1/(2\alpha+1)}$, (5.2) is satisfied if $\varepsilon_n$ satisfies $\sum_{i=1}^k \theta_{0,i}^2/\sigma_{i,k}^2 \lesssim n\varepsilon_n^2$ and $\sum_{i>k} \theta_{0,i}^2 \leq \varepsilon_n^2/2$.

It follows that the posterior concentrates at $\theta_0$ at the rate $\varepsilon_n$ that satisfies these requirements as well as the condition $\varepsilon_n \gtrsim n^{-\alpha/(2\alpha+1)}$. If the true parameter $\theta_0$ satisfies $\sum_{i=1}^\infty \theta_{0,i}^2 i^{2\alpha} < \infty$, then all three inequalities are satisfied for $\varepsilon_n$ a multiple of $n^{-\alpha/(2\alpha+1)}$. The rate $n^{-\alpha/(2\alpha+1)}$ is the minimax rate for this problem. □

Our prior is dependent on $n$, but with some more effort, it can be seen that the same conclusion can be obtained with a mixture prior of the form $\sum_n \lambda_n \Pi_n$ for suitable $\lambda_n$.

7.7. *Nonparametric regression with Gaussian errors.* Consider the nonparametric regression model, where we observe independent random variables $X_1, \ldots, X_n$ distributed as $X_i = f(z_i) + \varepsilon_i$ for an unknown regression function $f$, deterministic real-valued covariates $z_1, \ldots, z_n$ and normally distributed error variables $\varepsilon_1, \ldots, \varepsilon_n$ with zero means and variances $\sigma^2$. For simplicity, we assume that the error variance $\sigma^2$ is known. We also suppose that the covariates take values in a fixed compact set, which we will take as the unit interval, without loss of generality.

Let $f_0$ denote the true value of the regression function, let $P_{f,i}$ be the distribution of $X_i$ and let $P_f^{(n)}$ be the distribution of $(X_1, \ldots, X_n)$. Thus, $P_{f,i}$ is the normal measure with mean $f(z_i)$ and variance $\sigma^2$. Let $\mathbb{P}_n^z = n^{-1} \sum_{i=1}^n \delta_{z_i}$ be the empirical measure of the covariates and let $\|\cdot\|_n$ denote the norm on $L_2(\mathbb{P}_n^z)$.

By easy calculations, $K(P_{f_0,i}, P_{f,i}) = |f_0(z_i) - f(z_i)|^2/(2\sigma^2)$ and $V_{2,0}(P_{f_0,i}, P_{f,i}) = |f_0(z_i) - f(z_i)|^2/\sigma^2$ for all $i = 1, 2, \ldots, n$, whence the average Kullback–Leibler divergence and variance are bounded by a multiple of $\|f_0 - f\|_n^2/\sigma^2$ and hence it is enough to quantify prior concentration in $\|\cdot\|_n$-balls. The average Hellinger distance, as used in Theorem 4, is bounded above by $\|\cdot\|_n$, but is equivalent to this norm only if the class of regression functions is uniformly bounded, which makes it less attractive. However, it can be verified (cf. [5]) that the likelihood ratio test for $f_0$ versus $f_1$ satisfies the conclusion of Lemma 2 relative to $\|\cdot\|_n$ (instead of $d_n$ and $\theta_i = f_i$). Therefore, we may use the norm $\|\cdot\|_n$ instead of the average Hellinger distance throughout.

We shall construct priors based on series representations that are appropriate if $f_0 \in C^\alpha[0,1]$, where $\alpha > 0$ could be fractional. This means that $f_0$ is $\alpha_0$ times continuously differentiable with $\|f_0\|_\alpha < \infty$, $\alpha_0$ being the greatest integer less than $\alpha$ and the seminorm being defined by

$$\|f\|_\alpha = \sup_{x \neq x'} \frac{|f^{(\alpha_0)}(x) - f^{(\alpha_0)}(x')|}{|x - x'|^{\alpha - \alpha_0}}. \tag{7.9}$$

7.7.1. *Splines.* Fix an integer $q$ with $q \geq \alpha$. For a given natural number $K$, which will increase with $n$, partition the interval $(0,1]$ into $K$ subintervals $((k-1)/K, k/K]$ for $k = 1, 2, \ldots, K$. The space of splines of order $q$ relative to this partition is the collection of all functions $f : (0,1] \to \mathbb{R}$ that are $q - 2$ times continuously differentiable throughout $(0,1]$ and, if restricted to a subinterval $((k-1)/K, k/K]$, are polynomials of degree strictly less than $q$. These splines form a $J = (q + K - 1)$-dimensional linear space, with a convenient basis $B_1, B_2, \ldots, B_J$ being the B-splines, as defined in, for example, [11]. The B-splines satisfy (i) $B_j \geq 0$, $j = 1, 2, \ldots, J$, (ii) $\sum_{j=1}^{J} B_j = 1$, (iii) $B_j$ is supported inside an interval of length $q/K$ and (iv) at most $q$ of $B_1(x), \ldots, B_J(x)$ are nonzero at any given $x$. Let $B(z) = (B_1(z), \ldots, B_J(z))^T$ and write $\beta^T B$ for the function $z \mapsto \sum_j \beta_j B_j(z)$.

The basic approximation property of splines proved in [11], page 170, shows that for some $\beta_\infty \in \mathbb{R}^J$ (dependent on $J$),

$$\|\beta_\infty^T B - f_0\|_\infty \lesssim J^{-\alpha} \|f_0\|_\alpha. \tag{7.10}$$

Thus, by increasing $J$ appropriately with the sample size, we may view the space of splines as a sieve for the construction of the maximum likelihood estimator, as in Stone [28, 29], and for Bayes estimates as in [14, 15] for the problem of density estimation.

To put a prior on $f$, we represent it as $f_\beta(z) = \beta^T B(z)$ and induce a prior on $f$ from a prior on $\beta$. Ghosal, Ghosh and van der Vaart [14], in the context of density estimation, choose $\beta_1, \ldots, \beta_J$ i.i.d. uniform on an interval $[-M, M]$, the restriction to a finite interval being necessary to avoid densities with arbitrarily small values. In the present regression situation, a restriction to a compact interval is unnecessary and we shall choose $\beta_1, \ldots, \beta_J$ to be a sample from the standard normal distribution.

We need the regressors $z_1, z_2, \ldots, z_n$ to be sufficiently regularly distributed in the interval $[0, 1]$. In view of the spatial separation property of the B-spline functions, the precise condition can be expressed in terms of the covariance matrix $\Sigma_n = (\int B_i B_j \, d\mathbb{P}_n^z)$, namely

$$J^{-1}\|\beta\|^2 \lesssim \beta^T \Sigma_n \beta \lesssim J^{-1}\|\beta\|^2, \tag{7.11}$$

where $\|\cdot\|$ is the Euclidean norm on $\mathbb{R}^J$.

Under condition (7.11), we have that for all $\beta_1, \beta_2 \in \mathbb{R}^J$,

$$C\|\beta_1 - \beta_2\| \leq \sqrt{J}\|f_{\beta_1} - f_{\beta_2}\|_n \leq C'\|\beta_1 - \beta_2\| \tag{7.12}$$

for some constants $C$ and $C'$. This enables us to perform all calculations in terms of the Euclidean norms on the spline coefficients.

THEOREM 12. *Assume that the true density $f_0$ satisfies (7.10) for some $\alpha \geq \frac{1}{2}$, let (7.11) hold and let $\Pi_n$ be priors induced by a $N_J(0, I)$ distribution on the spline coefficients. If $J = J_n \sim n^{1/(1+2\alpha)}$, then the posterior converges at the minimax rate $n^{-\alpha/(1+2\alpha)}$ relative to $\|\cdot\|_n$.*

PROOF. We verify the conditions of Theorem 4. Let $f_{\beta_n}$ be the $L_2(\mathbb{P}_n^z)$-projection of $f_0$ onto the $J$-dimensional space of splines $f_\beta = \beta^T B$. Then $\|f_{\beta_n} - f_\beta\|_n \leq \|f_0 - f_\beta\|_n$ for every $\beta \in \mathbb{R}^J$ and hence, by (7.12), for every $\varepsilon > 0$, we have $\{\beta : \|f_\beta - f_0\|_n \leq \varepsilon\} \subset \{\beta : \|\beta - \beta_n\| \leq C'\sqrt{J}\varepsilon\}$. It follows that the $\varepsilon$-covering numbers of the set $\{f_\beta : \|f_\beta - f_0\|_n \leq \varepsilon\}$ for $\|\cdot\|_n$ are bounded by the $C\sqrt{J}\varepsilon$-covering numbers of a Euclidean ball of radius $C'\sqrt{J}\varepsilon$, which are of the order $D^J$ for some constant $D$. Thus, the entropy condition (3.2) is satisfied, provided that $J \lesssim n\varepsilon_n^2$.

By the projection property, with $\beta_\infty$ as in (7.10),

(7.13) $\qquad \|f_{\beta_n} - f_0\|_n \leq \|f_{\beta_\infty} - f_0\|_n \leq \|f_{\beta_\infty} - f_0\|_\infty \lesssim J^{-\alpha}.$

Combining this with (7.12) shows that there exists a constant $C''$ such that for every $\varepsilon \gtrsim 2J^{-\alpha}$, $\{\beta : \|f_\beta - f_0\|_n \leq \varepsilon\} \supset \{\beta : \|\beta - \beta_n\| \leq C''\sqrt{J}\varepsilon\}$. Together with the inclusion in the preceding paragraph and the definition of the prior, this implies that

$$\frac{\Pi_n(f : \|f - f_0\|_n \leq j\varepsilon)}{\Pi_n(f : \|f - f_0\|_n \leq \varepsilon)} \leq \frac{N_J(0, I)(\beta : \|\beta - \beta_n\| \leq C'j\sqrt{J}\varepsilon)}{N_J(0, I)(\beta : \|\beta - \beta_n\| \leq C''j\sqrt{J}\varepsilon)}$$

$$\leq \frac{N_J(0, I)(\beta : \|\beta\| \leq C'j\sqrt{J}\varepsilon)}{2^{-J/2}e^{-\|\beta_n\|^2} N_J(0, I)(\beta : \|\beta\| \leq C''j\sqrt{J}\varepsilon/\sqrt{2})}.$$

In the last step, we use Anderson's lemma to see that the numerator increases if we replace the centering $\beta_n$ by the origin, whereas to bound the denominator below, we use the fact that

$$\frac{dN_J(\beta_n, I)}{dN_J(0, I/2)}(\beta) = \frac{e^{-\|\beta - \beta_n\|^2/2}}{(\sqrt{2})^J e^{-\|\beta\|^2}} \geq 2^{-J/2} e^{-\|\beta_n\|^2}.$$

Here, by the triangle inequality, (7.12) and (7.13), we have that $\|\beta_n\| \lesssim \sqrt{J}\|f_{\beta_n}\|_n \lesssim \sqrt{J}(J^{-\alpha} + \|f_0\|_\infty) \lesssim \sqrt{J}$. Furthermore, the two Gaussian probabilities are left tail probabilities of the chi-square distribution with $J$ degrees of freedom. The quotient can be evaluated as

$$2^{J/2} e^{\|\beta_n\|^2} \frac{\int_0^{(C')^2 j^2 J\varepsilon^2} x^{J/2-1} e^{-x/2}\,dx}{\int_0^{(C'')^2 J\varepsilon^2/2} x^{J/2-1} e^{-x/2}\,dx}.$$

This is bounded above by $(Cj)^J$ for some constant $C$ if $\sqrt{J}\varepsilon$ remains bounded. Hence, to satisfy (3.4), it again suffices that $n\varepsilon_n^2 \gtrsim J$.

We conclude the proof by choosing $J = J_n \sim n^{1/(1+2\alpha)}$. $\square$

7.7.2. *Orthonormal series priors.* The arguments in the preceding subsection use the special nature of the B-spline basis only through the approximation inequality (7.10) and the comparison of norms (7.12). Theorem 12 thus extends to many other possible bases. One possibility is to use a sequence of orthonormal bases with good approximation properties for a given class of regression functions $f_0$. Then (7.11) should be replaced by

(7.14) $\qquad \|\beta_1 - \beta_2\| \lesssim \|f_{\beta_1} - f_{\beta_2}\|_n \lesssim \|\beta_1 - \beta_2\|.$

This is trivially true if the bases are orthonormal in $L_2(\mathbb{P}_n^z)$, but this requires that the basis functions change with the design points $z_1, \ldots, z_n$. One possible example is the discrete wavelet bases relative to the design points. All arguments remain valid in this setting.

7.8. *Binary regression.* Let $X_1, \ldots, X_n$ be independent observations with $P(X_i = 1) = 1 - P(X_i = 0) = F(\alpha + \beta z_i)$, where $z_i$ is a one-dimensional covariate, $\alpha$ and $\beta$ are parameters and $F$ is a cumulative distribution. Within the parametric framework, logit regression, where $F(z) = (1 + e^{-z})^{-1}$, or probit regression, where $F$ is the cumulative distribution function of the standard normal distribution, are usually considered. Recently, there has been interest in link functions of unknown functional form. The parameters $(F, \alpha, \beta)$ are separately not identifiable, unless some suitable restrictions on $F$ (such as given values of two quantiles of $F$) are imposed. For Bayesian estimation of $(F, \alpha, \beta)$, one therefore needs to put a prior on $F$ that conforms with the given restriction. However, in practice, one usually puts a Dirichlet process or a similar prior on $F$ and, independently of this, a prior on $(\alpha, \beta)$, and makes inference about, say, $z_0$, where $F(\alpha + \beta z_0) = 1/2$. Recently, Amewou-Atisso et al. [1] showed that the resulting posterior is consistent. In this section, we obtain the rate of convergence by an application of Theorem 4.

Because we directly measure distances between the distributions generating the data, identifiability issues need not concern us. The model and the prior can thus be described in a simpler form. We assume that $X_1, X_2, \ldots$ are independent Bernoulli variables, $X_i$ having success parameter $H(z_i)$ for an unknown, monotone link function $H$. As a prior on $H$, we use the Dirichlet process prior with base measure $\gamma((\cdot - \alpha)/\beta)$, for "hyperparameters" $(\alpha, \beta)$ distributed according to some given prior. This results in a mixture of Dirichlet process priors for $H$. Let the true value of $H$ be $H_0$, which is assumed to be continuous and nondecreasing.

In practice, $\gamma$ is often chosen to have support equal to the whole of $\mathbb{R}$ and $(\alpha, \beta)$ chosen to have support equal to $\mathbb{R} \times (0, \infty)$ so that the conditions on $\gamma$ and $(\alpha, \beta)$ described in the following theorem are satisfied.

THEOREM 13. *Assume that $z_1, z_2, \ldots, z_n$ lie in an interval $[a, b]$ strictly within the support of the true link function $H_0$ so that $H_0(a-) > 0$ and $H_0(b) < 1$. Let $H$ be the given mixture of Dirichlet process priors described previously with $\gamma$ and $(\alpha, \beta)$ having densities that are positive and continuous inside their supports. Assume that there exists a compact set $\mathbb{K}$ inside the support of the prior for $(\alpha, \beta)$ such that whenever $(\alpha, \beta) \in \mathbb{K}$, the support of the base measure $\gamma((\cdot - \alpha)/\beta)$ strictly contains the interval $[a, b]$. Then the posterior distribution of $H$ converges at the rate $n^{-1/3}(\log n)^{1/3}$ with respect to the distance $d_n$ given by (3.1).*

PROOF. Because the Hellinger distance between two Bernoulli distributions with success parameters $p$ and $q$ is equal to $(p^{1/2} - q^{1/2})^2 + ((1-p)^{1/2} - (1-q)^{1/2})^2$, we have

$$d_n^2(H_1, H_2) \leq \int |H_1^{1/2} - H_2^{1/2}|^2 \, d\mathbb{P}_n + \int |(1-H_1)^{1/2} - (1-H_2)^{1/2}|^2 \, d\mathbb{P}_n,$$

where $\mathbb{P}_n$ is the empirical distribution of $z_1, z_2, \ldots, z_n$. Both the classes $\{H^{1/2} : H \text{ is a c.d.f.}\}$ and $\{(1-H)^{1/2} : H \text{ is a c.d.f.}\}$ have $\varepsilon$-entropy bounded by a multiple of $\varepsilon^{-1}$, by Theorem 2.7.5 of [31]. Thus, any $\varepsilon_n \gtrsim n^{-1/3}$ satisfies (3.2).

By easy calculations, we have

$$K_i(H_0, H) = H_0(z_i) \log \frac{H_0(z_i)}{H(z_i)} + (1 - H_0(z_i)) \log \frac{1 - H_0(z_i)}{1 - H(z_i)},$$

$$V_{2;i}(H_0, H) \leq 2H_0(z_i)\left(\log \frac{H_0(z_i)}{H(z_i)}\right)^2 + 2(1 - H_0(z_i))\left(\log \frac{1 - H_0(z_i)}{1 - H(z_i)}\right)^2.$$

Under the conditions of the theorem, the numbers $H_0(z_i)$ are bounded away from 0 and 1. By Taylor's expansion, for any $\delta > 0$, there exists a constant $C$ (depending on $\delta$) such that

$$\sup_{\delta < p < 1-\delta} \sup_{q: |q-p| < \varepsilon} \left(p\left(\log \frac{p}{q}\right)^r + (1-p)\left(\log \frac{1-p}{1-q}\right)^r\right) \leq C\varepsilon^2, \qquad r = 1, 2.$$

Therefore, with $\|H - H_0\|_\infty = \sup\{|H(z) - H_0(z)| : z \in [a, b]\}$, we have $\max(n^{-1} \sum_{i=1}^n K_i(H_0, H), n^{-1} \sum_{i=1}^n V_{2;i}(H_0, H)) \lesssim \|H - H_0\|_\infty^2$. Hence, in order to satisfy (3.4), it suffices to lower bound the prior probability of the set $\{H : \|H - H_0\|_\infty \leq \varepsilon\}$.

For given $\alpha$ and $\beta$, the base measure is $\gamma((\cdot - \alpha)/\beta)$. For a given $\varepsilon > 0$, partition the line into $N \lesssim \varepsilon^{-1}$ intervals $E_1, E_2, \ldots, E_N$ such that $H_0(E_j) \leq \varepsilon$ and such that the $\gamma((\cdot - \alpha)/\beta)$-probability of every set $E_j$ (for $j = 1, 2, \ldots, N$) is between $A\varepsilon$ and 1 for a given positive constant $A$. Existence of such a partition follows from the continuity of $H_0$. It easily follows that for every $H$ such that $\sum_{j=1}^N |H(E_j) - H_0(E_j)| \leq \varepsilon$, we have $\|H - H_0\|_\infty \lesssim \varepsilon$. Furthermore, the conclusion is true even if $(\alpha, \beta)$ varies over $\mathbb{K}$. By Lemma 6.1 of [14], the prior probability of the set of all $H$ satisfying $\sum_{j=1}^N |H(E_j) - H_0(E_j)| \leq \varepsilon$ is at least $\exp(-c\varepsilon^{-1} \log \varepsilon^{-1})$ for some constant $c$. Furthermore, a uniform estimate works for all $(\alpha, \beta) \in \mathbb{K}$. Hence, (3.4) holds for $\varepsilon_n$, the solution of $n\varepsilon^2 = \varepsilon^{-1} \log \varepsilon^{-1}$, or for $\varepsilon_n = n^{-1/3}(\log n)^{1/3}$, which is only slightly weaker than the minimax rate $n^{-1/3}$. □

7.9. *Interval censoring.* Let $T_1, T_2, \ldots, T_n$ constitute an i.i.d. sample from a life distribution $F$ on $(0, \infty)$, which is subject to interval censoring by intervals $(l_1, u_1), \ldots, (l_n, u_n)$. We assume that the intervals are either nonstochastic or else we work conditionally on the realized values. Putting $(\delta_1, \eta_1), \ldots, (\delta_n, \eta_n)$, where $\delta_i = 1\{T_i \leq l_i\}$ and $\eta_i = 1\{l_i < T_i < u_i\}$, $i = 1, 2, \ldots, n$, the likelihood is given by $\prod_{i=1}^n (F(l_i))^{\delta_i} (F(u_i) - F(l_i))^{\eta_i} (1 - F(u_i))^{1-\delta_i-\eta_i}$. We may put the Dirichlet process prior on $F$. Under mild assumptions on the true $F_0$ and the base measure, the convergence rate under $d_n$ turns out to be $n^{-1/3}(\log n)^{1/3}$, which is the minimax rate, except for the logarithmic factor. Here, we use monotonicity of $F$ to bound the $\varepsilon$-entropy by a multiple of $\varepsilon^{-1}$ and we estimate prior probability concentration as $\exp(-c\varepsilon^{-1} \log \varepsilon^{-1})$ using methods similar to those used in the previous subsection. The details are omitted.

**8. Proofs.** In this section, we collect a number of technical proofs. For the proofs of the main results, we first present two lemmas.

LEMMA 9. *Let $d_n$ and $e_n$ be semimetrics on $\Theta$ for which tests satisfying the conditions of (2.2) exist. Suppose that for some nonincreasing function $\varepsilon \mapsto N(\varepsilon)$ and some $\varepsilon_n \geq 0$,*

$$(8.1) \quad N\left(\frac{\varepsilon \xi}{2}, \{\theta \in \Theta : d_n(\theta, \theta_0) < \varepsilon\}, e_n\right) \leq N(\varepsilon) \qquad \text{for all } \varepsilon > \varepsilon_n.$$

*Then for every $\varepsilon > \varepsilon_n$, there exist tests $\phi_n$, $n \geq 1$, (depending on $\varepsilon$) such that $P_{\theta_0}^{(n)} \phi_n \leq N(\varepsilon) \frac{e^{-Kn\varepsilon^2}}{1-e^{-Kn\varepsilon^2}}$ and $P_\theta^{(n)}(1 - \phi_n) \leq e^{-Kn\varepsilon^2 j^2}$ for all $\theta \in \Theta$ such that $d_n(\theta, \theta_0) > j\varepsilon$ and for every $j \in \mathbb{N}$.*

PROOF. For a given $j \in \mathbb{N}$, choose a maximal set of points in $\Theta_j = \{\theta \in \Theta : j\varepsilon < d_n(\theta, \theta_0) \leq (j+1)\varepsilon\}$ with the property that $e_n(\theta, \theta') > j\varepsilon\xi$ for every pair of points in the set. Because this set of points is a $j\varepsilon\xi$-net over $\Theta_j$ for $e_n$ and because $(j+1)\varepsilon \leq 2j\varepsilon$, this yields a set $\Theta'_j$ of at most $N(2j\varepsilon)$ points, each at $d_n$-distance at least $j\varepsilon$ from $\theta_0$, and every $\theta \in \Theta_j$ is within $e_n$-distance $j\varepsilon\xi$ of at least one of these points. (If $\Theta_j$ is empty, we take $\Theta'_j$ to be empty also.) By assumption, for every point $\theta_1 \in \Theta'_j$, there exists a test with the properties as in (2.2), but with $\varepsilon$ replaced by $j\varepsilon$. Let $\phi_n$ be the maximum of all tests attached in this way to some point $\theta_1 \in \Theta'_j$ for some $j \in \mathbb{N}$. Then

$$P_{\theta_0}^{(n)} \phi_n \leq \sum_{j=1}^\infty \sum_{\theta_1 \in \Theta'_j} e^{-Knj^2\varepsilon^2} \leq \sum_{j=1}^\infty N(2j\varepsilon) e^{-Knj^2\varepsilon^2} \leq N(\varepsilon) \frac{e^{-Kn\varepsilon^2}}{1 - e^{-Kn\varepsilon^2}}$$

and for every $j \in \mathbb{N}$,

$$\sup_{\theta \in \bigcup_{i>j} \Theta_i} P_\theta^{(n)}(1-\phi_n) \leq \sup_{i>j} e^{-Kni^2\varepsilon^2} \leq e^{-Knj^2\varepsilon^2},$$

where we have used the fact that for every $\theta \in \Theta_i$, there exists a test $\phi$ with $\phi_n \geq \phi$ and $P_\theta^{(n)}(1-\phi) \leq e^{-Kni^2\varepsilon^2}$. This concludes the proof. $\square$

LEMMA 10. *For $k \geq 2$, every $\varepsilon > 0$ and every probability measure $\bar{\Pi}_n$ supported on the set $B_n(\theta_0, \varepsilon; k)$, we have, for every $C > 0$,*

$$(8.2) \qquad P_{\theta_0}^{(n)}\left(\int \frac{p_\theta^{(n)}}{p_{\theta_0}^{(n)}}\, d\bar{\Pi}_n(\theta) \leq e^{-(1+C)n\varepsilon^2}\right) \leq \frac{1}{C^k(n\varepsilon^2)^{k/2}}.$$

PROOF. By Jensen's inequality applied to the logarithm, with $l_{n,\theta} = \log(p_\theta^{(n)}/p_{\theta_0}^{(n)})$, we have $\log \int (p_\theta^{(n)}/p_{\theta_0}^{(n)})\, d\bar{\Pi}_n(\theta) \geq \int l_{n,\theta}\, d\bar{\Pi}_n(\theta)$. Thus, the probability in (8.2) is bounded above by

$$(8.3) \qquad P_{\theta_0}^{(n)}\left(\int (l_{n,\theta} - P_{\theta_0}^{(n)} l_{n,\theta})\, d\bar{\Pi}_n(\theta) \leq -n(1+C)\varepsilon^2 - \int P_{\theta_0}^{(n)} l_{n,\theta}\, d\bar{\Pi}_n(\theta)\right).$$

For every $\theta \in B_n(\theta_0, \varepsilon; k)$, we have $P_{\theta_0}^{(n)} l_{n,\theta} = -K(p_{\theta_0}^{(n)}, p_\theta^{(n)}) \geq -n\varepsilon^2$. Consequently, by Fubini's theorem and the assumption that $\bar{\Pi}_n$ is supported on this set, the expression on the right-hand side of (8.3) is bounded above by $-Cn\varepsilon^2$. An application of Markov's inequality yields the upper bound

$$\frac{P_{\theta_0}^{(n)}|\int (l_{n,\theta} - P_{\theta_0}^{(n)} l_{n,\theta})\, d\bar{\Pi}_n(\theta) \wedge 0|^k}{(Cn\varepsilon^2)^k} \leq \frac{P_{\theta_0}^{(n)} \int |l_{n,\theta} - P_{\theta_0}^{(n)} l_{n,\theta}|^k\, d\bar{\Pi}_n(\theta)}{(Cn\varepsilon^2)^k},$$

by another application of Jensen's inequality. The right-hand side is bounded by $C^{-k}(n\varepsilon^2)^{-k/2}$, by the assumption on $\bar{\Pi}_n$. This concludes the proof. $\square$

PROOF OF THEOREM 1. By Lemma 9, applied with $N(\varepsilon) = \exp(n\varepsilon_n^2)$ (constant in $\varepsilon$) and $\varepsilon = M\varepsilon_n$ in its assertion, where $M \geq 2$ is a large constant to be chosen later, there exist tests $\phi_n$ that satisfy $P_{\theta_0}^{(n)}\phi_n \leq e^{n\varepsilon_n^2}(1 - e^{-KnM^2\varepsilon_n^2})^{-1}e^{-KnM^2\varepsilon_n^2}$ and $P_\theta^{(n)}(1-\phi_n) \leq e^{-KnM^2\varepsilon_n^2 j^2}$ for all $\theta \in \Theta_n$ such that $d_n(\theta, \theta_0) > M\varepsilon_n j$ and for every $j \in \mathbb{N}$. The first assertion implies that if $M$ is sufficiently large to ensure that $KM^2 - 1 > KM^2/2$, then as $n \to \infty$, for any $J \geq 1$, we have

$$(8.4) \qquad P_{\theta_0}^{(n)}[\Pi_n(d_n(\theta, \theta_0) \geq JM\varepsilon_n | X^{(n)})\phi_n] \leq P_{\theta_0}^{(n)}\phi_n \lesssim e^{-KM^2 n\varepsilon_n^2/2}.$$

Setting $\Theta_{n,j} = \{\theta \in \Theta_n : M\varepsilon_n j < d_n(\theta, \theta_0) \leq M\varepsilon_n(j+1)\}$ and using (2.2), we obtain, by Fubini's theorem,

$$(8.5) \qquad P_{\theta_0}^{(n)}\left[\int_{\Theta_{n,j}} \frac{p_\theta^{(n)}}{p_{\theta_0}^{(n)}}\, d\Pi_n(\theta)\, (1-\phi_n)\right] \leq e^{-KnM^2\varepsilon_n^2 j^2} \Pi_n(\Theta_{n,j}).$$

POSTERIOR CONVERGENCE RATES 31Fix some $C > 0$. By Lemma 10, we have, on an event $A_n$ with probability at least $1 - C^{-k}(n\varepsilon_n^2)^{-k/2}$,

$$\int \frac{p_\theta^{(n)}}{p_{\theta_0}^{(n)}} d\Pi_n(\theta) \geq \int_{B_n(\theta_0,\varepsilon_n;k)} \frac{p_\theta^{(n)}}{p_{\theta_0}^{(n)}} d\Pi_n(\theta) \geq e^{-(1+C)n\varepsilon_n^2} \Pi_n(B_n(\theta_0, \varepsilon_n; k)).$$

Hence, decomposing $\{\theta \in \Theta : d_n(\theta, \theta_0) > JM\varepsilon_n\} = \cup_{j \geq J} \Theta_{n,j}$ and using (8.5), the last display and (2.5), we have, for every sufficiently large $J$,

$$P_{\theta_0}^{(n)}[\Pi_n(\theta \in \Theta_n : d_n(\theta, \theta_0) > J\varepsilon_n M | X^{(n)})(1 - \phi_n) 1_{A_n}]$$
$$\leq \sum_{j \geq J} e^{-n\varepsilon_n^2 (KM^2 j^2 - 1 - C - \frac{1}{2}KM^2 j^2)},$$

by assumption (2.5). This converges to zero as $n \to \infty$ for fixed $C$ and fixed, sufficiently large $M$ and $J$ if $n\varepsilon_n^2 \to \infty$; it converges to zero for fixed $M$ and $C$ as $J = J_n \to \infty$ if $n\varepsilon_n^2$ is bounded away from zero.

Combining the preceding results, we have, for sufficiently large $M$ and $J$,

$$(8.6) \quad \begin{aligned} &P_{\theta_0}^{(n)} \Pi_n(\theta \in \Theta : d_n(\theta, \theta_0) > M\varepsilon_n J | X^{(n)}) \\ &\leq \frac{1}{C^k (n\varepsilon_n^2)^{k/2}} + 2e^{-KM^2 n\varepsilon_n^2 / 2} + \sum_{j \geq J} e^{-n\varepsilon_n^2 (\frac{1}{2}KM^2 j^2 - 1 - C)}. \end{aligned}$$

The rest of the conclusion follows easily; see the proof of Theorem 2.4 of [14]. □

PROOF OF THEOREM 2. If $\varepsilon_n \gtrsim n^{-\alpha}$ and $k(1 - 2\alpha) > 2$ for $\alpha \in (0, 1/2)$, then $n\varepsilon_n^2 \to \infty$ and $\sum_{n=1}^\infty (n\varepsilon_n^2)^{-k/2} < \infty$. For $C = 1/2$, the first term on the right-hand side of (8.6) dominates and the sum over $n$ of the terms in (8.6) converges. The result (i) follows by the Borel–Cantelli lemma.

For assertion (ii), we note that $\varepsilon_n \gtrsim n^{-\alpha}$ and $k(1 - 2\alpha) \geq 4\alpha$ together imply that $(n\varepsilon_n^2)^{-k/2} \lesssim \varepsilon_n^2$. The other terms are exponentially small. □

PROOF OF LEMMA 1. Because $P_{\theta_0}^{(n)}(p_\theta^{(n)}/p_{\theta_0}^{(n)}) \leq 1$, Fubini's theorem implies that $P_{\theta_0}^{(n)}[\int_{\Theta \setminus \Theta_n} (p_\theta^{(n)}/p_{\theta_0}^{(n)}) d\Pi_n(\theta)] \leq \Pi_n(\Theta \setminus \Theta_n)$. Let the events $A_n$ be as in the proof of Theorem 1, so that the denominator of the posterior is bounded below by $e^{-(1+C)n\varepsilon_n^2} \Pi_n(B_n(\theta_0, \varepsilon_n; k))$ on $A_n$. Combining this with the preceding display gives

$$P_{\theta_0}^{(n)}[\Pi_n(\theta \notin \Theta_n | X^{(n)}) 1_{A_n}] \leq \frac{\Pi_n(\Theta \setminus \Theta_n) e^{(1+C)n\varepsilon_n^2}}{\Pi_n(B_n(\theta_0, \varepsilon_n; k))} \leq o(1) e^{-n\varepsilon_n^2(1-C)},$$

by the assumption on $\Pi_n(\Theta \setminus \Theta_n)$. The rest of the proof can be completed along the lines of that of Theorem 2.4 of [14]. □

Department of Statistics
North Carolina State University
2501 Founders Drive
Raleigh, North Carolina 27695-8203
USA
E-mail: ghosal@stat.ncsu.edu

Division of Mathematics
and Computer Science
Vrije Universiteit
De Boelelaan 1081a
1081 HV Amsterdam
The Netherlands
E-mail: aad@cs.vu.nl